\documentclass[10pt,journal]{IEEEtran}
\usepackage{amsmath}
\usepackage{amsfonts}
\usepackage{amssymb}
\usepackage{graphicx}
\DeclareGraphicsExtensions{.pdf,.png,.jpg}
\usepackage{algorithmic}
\usepackage{bm}
\usepackage{xcolor} 
\usepackage{graphicx,epstopdf,tikz}
 \DeclareMathAlphabet{\itbf}{OML}{cmm}{b}{it}
\DeclareMathAlphabet\mathbfcal{OMS}{cmsy}{b}{n}

\renewcommand{\tilde}{\widetilde}

\def\RR{\mathbb{R}}
\def\bx{{{\itbf x}}}
\def\bz{{{\itbf z}}}
\def\by{{{\itbf y}}}

\def\bphi{{\boldsymbol{\varphi}}}

\def\balpha{{\boldsymbol{\alpha}}}

\def\eps{\varepsilon}

\def\bR{{\itbf R}}

\def\bV{{\itbf V}}
\def\vE{\vec{\itbf E}}
\def\vH{\vec{\itbf H}}
\def\bH{{\itbf H}}
\def\bU{{\itbf U}}

\def\bI{{\itbf I}}

\def\cL{{\mathcal{L}}}

\def\m0{\underline{\underline{\bf 0}}}

\def\lb{\left <}

\def\rb{\right >}

\def\cF{\mathcal{F}}
\def\cT{\mathcal{T}}
\def\om{\omega}
\def\la{\lambda}
\def\12{{\frac{1}{2}}}

\usepackage[normalem]{ulem}
\usepackage{pgfplots}
\usepackage{tikz}
\newlength{\figurewidth}
\newlength{\figureheight}

\newtheorem{thm}{Theorem}

\usepackage{calrsfs} 
\newcommand{\FIGDIR}{FiguresPaper}
\newcommand{\IMGDIR}{FiguresPaper}

\begin{document}

\title{Quantitative synthetic aperture radar inversion}
\author{Liliana Borcea, Josselin Garnier, Alexander V. Mamonov, J\"{o}rn Zimmerling

\thanks{Borcea is with the Applied Physics \& Applied Math. department, Columbia University, New York, NY 10027 USA (e-mail: lb3539@columbia.edu). }
\thanks{Garnier is with CMAP, CNRS, Ecole polytechnique, Institut Polytechnique de Paris, 91120 Palaiseau, France (e-mail: 
josselin.garnier@polytechnique.edu)}
\thanks{Mamonov is with the Mathematics department, University of Houston, Houston, TX 77204-3008 USA (e-mail: 
mamonov@math.uh.edu)}
\thanks{Zimmerling is with the department of Information Technology, Division of Scientific Computing, Uppsala Universitet, 75105 Uppsala,
Sweden (e-mail: jorn.zimmerling@it.uu.se)}
}

\markboth{submitted to: IEEE Transactions on Antennas and Propagation}
{Shell \MakeLowercase{\textit{Borcea, Garnier, Mamonov, Zimmerling}}: Inverse scattering for SAR}
\maketitle

\begin{abstract} We study an inverse scattering problem for monostatic synthetic aperture radar (SAR): Estimate the wave speed in a heterogeneous, isotropic and nonmagnetic medium probed by waves emitted and measured by a moving antenna.  The forward map, from the wave speed to the measurements, is   derived from Maxwell's equations. It is a nonlinear map that accounts for multiple scattering and  it is very oscillatory at high frequencies. This makes the standard, nonlinear least squares data fitting formulation of the inverse problem difficult to solve. We introduce an alternative, two-step approach: The first step computes the nonlinear map from the measurements to an approximation of the electric field  inside the unknown medium aka, the internal wave. 
This is done for each antenna location in a non-iterative manner.  
%This is done for each antenna location, in an efficient, non-iterative manner.  
The  internal wave fits the data by construction, but it does not solve Maxwell's equations.  The second step uses optimization to minimize the discrepancy between the internal  wave  and the solution of Maxwell's equations,  for all antenna locations.  The optimization is iterative.
The first step defines an imaging function whose computational cost is comparable to that of standard SAR imaging, but it  gives a better estimate of the support of targets. Further iterations improve the quantitative estimation of the wave speed. We assess the performance of the method with numerical simulations and compare the results with those of standard inversion.

\end{abstract}

% Note that keywords are not normally used for peer review papers.

\section{Introduction}

In monostatic synthetic aperture radar (SAR), a moving  antenna  probes an unknown heterogeneous medium by emitting pulsed or chirped, 
directed beams of electromagnetic radiation and then measuring the generated electric field. The measurements are  processed  with methods like matched filtering to obtain an imaging function, aka an image \cite{Cheney,CheneyB,Curlander}.  This function is evaluated at  points in the imaging domain and it gives a qualitative description of the medium: 
%Ideally,
It peaks near the supports of  reflectors (targets),  but it does not quantify the properties of the medium, modeled in Maxwell's equations by spatially variable and unknown coefficients. 

We are interested in the quantitative estimation of isotropic and nonmagnetic media, modeled by the  positive dielectric permitivity function $\eps(\bx)$ and  the constant magnetic permeability $\mu$. These define the wave speed $c(\bx) = 1/\sqrt{\mu \eps(\bx)}$. The SAR data $\{D_s(t), t \in \cT, ~ 1\le s \le  S\}$ depend on two variables: The slow time index $s$, which counts the antenna location on the flight track at the time of emission, and the fast time $t $ that runs between emissions, in the time interval $\cT \subset \RR$. It is impossible to obtain three-dimensional estimates of $\eps(\bx)$  from such data.  Therefore, the inversion is restricted 
to two dimensions by either assuming that $\eps(\bx)$ does not change in one direction, or by inverting on a known surface \cite{Cheney,Gilman}.

The estimation of $\eps(\bx)$ requires inverting in some appropriate sense the nonlinear forward map
\begin{equation}
\eps(\bx) \stackrel{\cF}{\mapsto} \{D_s(t),  ~t \in \cT, 1 \le s \le S\},
\label{eq:1}
\end{equation}
defined by evaluating the solution of  Maxwell's equations at the successive antenna locations.
 All SAR imaging is based on the single-scattering approximation of 
\eqref{eq:1}, which linearizes $\cF(\eps)$ about the constant and known permittivity $\eps_o$ near the antenna  \cite{Cheney,CheneyB,Gilman}. The linearization is convenient to use, but it is inadequate for describing the wave field in the presence of extended targets \cite{Gilman,Symes,Marks}, because it models them as a superposition of point scatterers that do not interact with each other. 
Consequently, SAR images can contain artifacts that complicate tasks like automatic target recognition \cite{TargetRec}.  

There are various improvements of the single-scattering approximation, that are still  based on weak scattering and scale separation assumptions \cite{Symes,Gilman,Marks}.
However, to obtain  a quantitative estimate 
of $\eps(\bx)$, one needs to take multiple scattering into account and ``invert" the nonlinear map $\cF$. Since $\cF$ is never invertible in the strict mathematical sense, the typical approach is to solve the nonlinear least squares data fitting optimization over the search permittivity $\hat{\eps}(\bx)$,
\begin{equation}
\min_{\hat{\eps}} \sum_{s=1}^{S} \int_{\cT} \, dt \left| \big[\cF(\hat{\eps})\big]_s(t) - D_s(t)\right|^2.
\label{eq:2}
\end{equation}
This formulation is called full waveform inversion (FWI) in the geophysics community \cite{Virieux}, but it is relevant elsewhere, including SAR. The optimization~\eqref{eq:2}  is augmented with some regularization penalty on $\hat{\eps}(\bx)$ and it is solved with a gradient based  iterative approach, where $\cF$ is linearized at each step. In particular, if one  starts with  $\hat{\eps}(\bx) = \eps_o$, the linearization is  the same as the Born approximation used in SAR imaging and  the 
first iterate is the solution of the normal equation for linear least squares data fitting.

High-frequency asymptotic analysis \cite{Beylkin,Symes} shows that the Born data model can be written approximately as a generalized Radon transform of the  ``reflectivity" function
{$\eps(\bx)/\eps_o-1$}
 that can be inverted via filtered backprojection \cite{Beylkin,Nolan}. 
Moreover, the normal operator behaves microlocally as an  identity operator \cite{Beylkin,Rakesh}. This holds approximately, if the aperture and bandwidth are large enough. It means that by applying the adjoint of the linearized forward map to the data, we obtain an image that is qualitatively the same as the solution of the normal equation. The adjoint map is known as a matched filter \cite{CheneyB} and the approach, 
called matched filtering,  is commonly used in SAR imaging and elsewhere. However, it does not produce a quantitative estimate 
of $\eps(\bx)$ and, as mentioned above, the images can have artifacts due to the neglected multiple scattering effects that are present in the data 
but are not accounted for in imaging.  

Iterations for \eqref{eq:2} may not give a good result, 
because the optimization can get stuck  in one of the many spurious minima of the objective function, that arise 
far and near the true $\eps(\bx)$ \cite{Virieux}.  This happens especially at the high frequencies used  in SAR, 
because the forward map is quite oscillatory \cite{Barucq}. The bad behavior of the data fitting objective function 
has motivated much research into alternative formulations, like using a different norm to quantify the 
data misfit \cite{Engquist1,Engquist2}, expanding the search space \cite{Huang,Herrman} and boundary control \cite{Belishev}. 
The first two approaches have shown some success in seismic imaging, but have not been tried for SAR data. Boundary control 
requires different and usually unavailable measurements i.e., the Dirichlet to Neumann map.   

In this paper we propose a different approach to the quantitative estimation of $\eps(\bx)$. It is rooted in the recent
data driven reduced order model (ROM) methods  for inverse scattering with multiple input multiple output (MIMO) active arrays \cite{ROM1}-\cite{ROM5}.
The ROMs are physics and data driven algebraic models (matrices) of the wave operator \cite{ROM3} or the propagator operator that maps snapshots of the wave field from one time instant to the next \cite{ROM1,ROM3,ROM4}. They were originally defined for  the second-order acoustic wave equation and were used to determine either $c(\bx)$, assuming a constant density \cite{ROM2,ROM3}, or determine the density assuming that $c(\bx)$ is constant  \cite{ROM1,ROM6,Druskin1}. ROMs for first-order hyperbolic systems with 
multiple unknown coefficients were introduced recently in \cite{ROM5}.  An important idea that emerged from these studies is that  the ROM propagator  can be used to approximate the wave field inside the medium aka, the internal wave \cite{ROM2,Druskin3}.
This was used in the Lippmann-Schwinger integral equation for the scattered wave field to linearize approximately the data fitting inversion process, 
while still taking multiple scattering into account \cite{Druskin2,ROM2}. It was also used in \cite{Druskin1} to map monostatic SAR data to fixed array MIMO data. These procedures work well under two conditions: (1)
The kinematics of the medium (the smooth part of $c(\bx)$) is not strongly perturbed.  For example, $c(\bx) = c_o = 1/\sqrt{\mu \eps_o}$  in \cite{Druskin1}. (2) The frequency of the probing signals is not too high.

It was observed in  \cite{ROM3} that the approximated internal wave is guaranteed to fit the MIMO data used to compute the ROM propagator. However, the 
{internal}
wave does not solve  the wave equation. This motivated an iterative optimization formulation 
of the inverse problem, that minimizes the misfit between the  approximated internal wave and the solution of the wave equation at the search speed  \cite{ROM3,ROM5}. In this paper, we extend the results in \cite{ROM3,ROM5} to monostatic SAR. The idea is to construct a data driven approximation of the internal electric  wave field for each slow time. The construction starts with Maxwell's equations excited by the antenna that emits a pulsed beam of electromagnetic radiation. The same antenna measures the generated electric field. We show how to obtain from the measurements an approximate internal wave, by adapting and improving the procedure in  \cite{ROM5}. Then, we couple the results for all the slow times, via optimization. 

The paper is organized as follows: We begin in section \ref{sect:setup} with the mathematical formulation of the inverse problem. The SAR data driven approximation of the internal wave and its properties are described  in section \ref{sect:inverse}. In section \ref{sect:numerics} 
we formulate the inverse problem as an optimization and describe the computational setup. The numerical results are in section \ref{sect:numerics_res}. We end with a summary in section \ref{sect:summary}.

\section{Formulation of the inverse problem}
\label{sect:setup}
Assume a
% two-dimensional 
{three-dimensional}
setting, where the medium is  invariant in the $z$ direction, orthogonal to the inversion plane with coordinates $\bx = (x,y)$. Consider  $E$-polarized  waves modeled by the electric  field
$
\vE(\bx,t) = \vec{\bz} E(\bx,t)$ and the magnetic field $
\vH(\bx,t) = \vec{{\itbf x}} H_x(\bx,t) + \vec{\itbf y} H_y(\bx,t),$
where $\vec{{\itbf x}}, \vec{{\itbf y}}$ and $\vec{{\itbf z}}$ are the unit vectors along the axes of the orthogonal coordinate system 
$(x,y,z)$. 

We index by the slow time $s$ the current-density forcing 
$J_s(\bx,t)$  from the antenna and the generated electromagnetic fields. These
satisfy the  first-order hyperbolic system 
\begin{equation}
\begin{pmatrix} 
\eps(\bx) \partial_t & \partial_y & - \partial_x \\
\partial_y & \mu \partial_t & 0 \\
- \partial_x & 0 & \mu \partial_t \end{pmatrix} 
\begin{pmatrix} 
E_s(\bx,t) \\ H_{x,s}(\bx,t) \\ H_{y,s}(\bx,t)
\end{pmatrix} = 
\begin{pmatrix} {J_s(\bx,t)}/{c_o} \\ 0 \\ 0 \end{pmatrix}, 
\label{eq:S3}
\end{equation}
derived from Maxwell's equations, where the normalization {of the source}
by  $c_o$  is used for convenience.
There is no wave prior to the excitation, so we set the fields to zero 
at time $t$ preceding the temporal support of $J_s(\bx,t)$.   

The antenna is modeled as a phased array that can  emit a probing beam. 
The direction and origin of the beam changes with the slow time. Thus,  we write $J_s(\bx,t)$ in the local coordinate system 
 $\bx_s = (x_s,y_s)$ with origin at the center of the antenna, and  with $y_s$ along the axis of the beam, rotated from $y$ by some angle $\theta_s$.  The relation $\bx = \bx(\bx_s)$ depends on the antenna flight path and the excitation is modeled by 
\begin{equation}
J_s(\bx(\bx_s),t) = \delta(y_s) b(x_s,t), \quad \bx_s = (x_s,y_s).
\label{eq:S2}
\end{equation}
Here $b(x_s,t)$ accounts for the cross-range profile of the beam at its origin and the emitted signal. For simplicity, we assume a pulse signal, supported at $ t \in [-T_b,T_b]$, modulated at the central frequency $\om_o$ and with bandwidth of order $1/T_b$.  
However, chirped signals can be accommodated, as well.

\vspace{-0.1in} \subsection{Symmetrization} 
We will use functional calculus on the wave operator, so it is convenient to rewrite 
the system~\eqref{eq:S3} in the form 
\begin{equation}
(\partial_t + \cL) \begin{pmatrix} u_s(\bx,t) \\ {\bH}_s(\bx,t) \end{pmatrix} = \begin{pmatrix} {J_s(\bx,t)} \\ {\bf 0}\end{pmatrix},
\label{eq:S4}
\end{equation}
where 
\begin{equation}
u_s(\bx,t) = E_s(\bx,t)/[\mu c(\bx)],\quad c(\bx) = \frac{1}{\sqrt{\mu \eps(\bx)}},
\label{eq:S4.1}
\end{equation}
 and $\bH_s(\bx,t)$ is the two-dimensional field with components $H_{x,s}(\bx,t)$ and 
 {$H_{y,s}(\bx,t)$}.
The operator 
\begin{equation}
\cL = \begin{pmatrix} 
0& c(\bx) \partial_y & - c(\bx) \partial_x \\
\partial_y[c(\bx) \cdot ] & 0& 0 \\
- \partial_x[c(\bx) \cdot] & 0 & 0\end{pmatrix} , 
\label{eq:S5}
\end{equation}
is skew-adjoint when acting on the space of sufficiently regular functions \cite{Monk}  with compact support in  $\RR^2$. Such functions are sufficient for our study, because the waves propagate at finite speed, so during the time  interval $\cT$ they are supported inside  some bounded ball in $\RR^2$, centered at the antenna.

\vspace{-0.1in} \subsection{Measurements and the inverse problem}
The measurements are modeled by 
\begin{equation}
D_s(t) = \int_{\mathbb{R}^2} d \bx \int_{-T_b}^{T_b} dt'  J_s(\bx,-t')  u_s(\bx,t-t'), 
\label{eq:S6}
\end{equation}
for $t \in \cT$ and $1 \le s \le S$, where $J_s(\bx,t)$ is supported at the antenna, per definition~\eqref{eq:S2}. Recall that $\eps(\bx) = \eps_o$ near the antenna. Thus, aside from the known constant $\sqrt{\mu/\eps_o}$,  equation~\eqref{eq:S6} gives 
the net electric field at the phased array modeling the antenna, convolved with the time-reversed probing signal. Such a convolution is commonly used  in 
SAR to compress long, chirped signals \cite{PulseComp}. 

\vskip 0.05in \noindent \textbf{Inverse problem:} Estimate from the data~\eqref{eq:S6} the permittivity $\eps(\bx)$ and therefore the speed  $c(\bx)$  in the compact domain  $\Omega_{\rm im} \subset \RR^2$, assumed to contain  the support of $\eps(\bx)-\eps_o$.

\section{The approximated internal wave}
\label{sect:inverse}
The inversion  uses a data driven approximation of the snapshots of $u_s(t,\bx)$. We  describe it  in this section.  

\vspace{-0.1in}
 \subsection{Snapshots of the wave field}
 \label{sect:defSnap}
The snapshots are defined on a uniform time grid 
\begin{equation}
m \tau + T, \quad  0 \le m \le M-1,
\label{eq:I1}
\end{equation}
with origin at $T  \gg  T_b$, 
%where $T$ is the time of travel of the pulsed beam  from the antenna to the imaging region $\Omega_{\rm im}$:
where $T$ is the time of travel of the pulsed beam  from the antenna to the imaging region $\Omega_{\rm im}$: $T= \inf_{\bx \in \Omega_{\rm im} ,\by \in {\rm supp}(J_s)} \| \bx- \by\| / c_o -T_b$.
The steps $\tau$ satisfy the Nyquist sampling criterium for the highest frequency in the bandwidth of $J_s(\bx,t)$ and 
\begin{equation}
M = \max\{m \in \mathbb{N}  ~\mbox{such that}~  2[T + (m-1) \tau] \in \cT\}.
\label{eq:I2}
\end{equation}
We consider the instances~\eqref{eq:I1}  because at $ t \le T$, the wave fields are not affected by the heterogeneity of the medium (supported in $\Omega_{\rm im}$) and are  thus  the same as those computed with permittivity $\eps_o$.

The ``primary" wave snapshots are defined by 
\begin{equation}
\bphi_{s,m}(\bx) = \begin{pmatrix} u_s(\bx,m \tau+T) \\ {\bH}_s(\bx,m \tau+T) \end{pmatrix} = e^{-m \tau \cL } \bphi_{s,0}(\bx),
\label{eq:I3}
\end{equation}
and their initial state is given by  the solution of~\eqref{eq:S4} at $t=T$
\begin{equation}
\bphi_{s,0}(\bx) = \int_{-T_b}^{T_b} dt' e^{-(T-t')\cL}\begin{pmatrix} {J_s(\bx,t')} \\ {\bf 0}\end{pmatrix}.
\label{eq:I4}
\end{equation}
Here we used functional calculus on the skew-adjoint operator $\cL$ and introduced the unitary evolution operator $e^{-t \cL}$.

To separate the scaled electric field (per~\eqref{eq:S4.1})  from the magnetic field, 
we introduce the ``adjoint" snapshots, defined via the ``time reversal" multiplication operator
\begin{equation}
 \mathbb{T}= \begin{pmatrix} 1 & {\bf 0}^T \\ 
 {
 {\bf 0}
 }
  & - \bI_2 \end{pmatrix},
 \label{eq:I5}
\end{equation}
where $\bI_2$ is the $2\times 2$ identity matrix. It is easy to check the commutation relation 
$
\cL \mathbb{T} = - \mathbb{T} \cL,
$
which implies that 
\begin{equation}
\mathbb{T} e^{-t \cL} = e^{t \cL} \mathbb{T}.
 \label{eq:I7}
 \end{equation}
The adjoint vectorial wave is 
\begin{equation}
\bphi^\star_{s,m}(\bx) = \mathbb{T} \bphi_{s,m}(\bx) \stackrel{\eqref{eq:I3},\eqref{eq:I7}}= e^{m \tau \cL } \bphi^\star_{s,0}(\bx).
\label{eq:I8}
\end{equation}
It evolves according to the  unitary operator $e^{t \cL}$, corresponding to the adjoint $-\cL$ of $\cL$, and its initial state is 
\begin{align}
\hspace{-0.1in}\bphi^\star_{s,0}(\bx) &= 
\mathbb{T} \bphi_{s,0}(\bx)  =  \int_{-T_b}^{T_b} dt' e^{(T-t')\cL}\begin{pmatrix} {J_s(\bx,t')} \\ {\bf 0}\end{pmatrix}.
\label{eq:I4_s}
\end{align}

\vspace{-0.1in}
 \subsection{From the data to inner products of the snapshots}
 \label{sect:approxSnap}
We wish to map the data ~\eqref{eq:S6} to an approximation of 
\begin{equation}
u_{s,m}(\bx)= u_s(\bx,m \tau + T), \quad 0 \le m \le M-1.
\label{eq:snapU}
\end{equation}
These are related to the primary snapshots~\eqref{eq:I3} and the adjoint snapshots~\eqref{eq:I8} by 
\begin{equation}
\begin{pmatrix} u_{s,m}(\bx) \\ {\bf 0} \end{pmatrix} =
\frac{1}{2} \left[ \bphi_{s,m}(\bx)+\bphi^\star_{s,m}(\bx)\right]. 
\label{eq:I9}
\end{equation}

The next theorem, proved in  appendix \ref{ap:A}, states that we  can determine from~\eqref{eq:S6}
the inner products of the snapshots, without knowing the medium. This allows us to define  in section \ref{sect:approxU} 
a family of 
approximations of~\eqref{eq:snapU} that all fit the measurements~\eqref{eq:S6} but are not necessarily solutions of the 
the hyperbolic system~\eqref{eq:S4}. To drive the approximation to the true snapshots and consequently, 
determine $\eps(\bx)$, we formulate the inversion as an optimization problem in section \ref{sect:numerics}.

\vskip 0.05in\begin{thm}
\label{prop.2}
\emph{Assume that the probing signal is even in time. Let $\lb \cdot, \cdot \rb$ denote the %Euclidean
  inner product 
\begin{equation}
\lb a, b \rb = \int_{\RR^2} a(\bx) \overline{b}(\bx) d\bx ,
\end{equation} and denote by 
$\mathbb{G}_s$ the $M \times M$  Gramian matrix with entries 
\begin{equation}
(\mathbb{G}_s)_{m,j} = \lb u_{s,m},u_{s,j}\rb, \quad 0 \le m,j \le M-1.
\label{eq:I10}
\end{equation}
 This is a symmetric matrix with Toeplitz plus Hankel structure.  Its entries above the diagonal are 
\begin{align}
(\mathbb{G}_s)_{m,m+j} &= \frac{1}{2} \left[ D_s(j\tau) + D_s(-j\tau)\right] \nonumber  \\
&+ \frac{1}{2}D_s(2 T + (2m+j)\tau),
\label{eq:I11}
\end{align}
for $0 \le m \le M-1$ and $0 \le j \le M-1-m$. The entries below the diagonal are obtained from symmetry.}
\end{thm}

\vskip 0.05in 
The assumption that the probing signal is even in time is used in Theorem \ref{prop.2} to simplify the formula~\eqref{eq:I11}.
If the signal is not even,  $D_s$ in the last term of ~\eqref{eq:I11} should be replaced by the analogue of~\eqref{eq:S6}, where the wave is convolved with $J_s(\bx,t)$, without time reversal (see appendix \ref{ap:A}).

Note that all the terms in ~\eqref{eq:I11} are measured, except for $D_s(-j\tau)$. The homogeneous initial condition 
 at $t < -T_b$ and the definition~\eqref{eq:S6} imply that 
$
D_s(-j \tau) \ne 0$ if  $0 \le j  \tau \le 2 T_b.$
The wave does not interact with the heterogeneity in the medium until $t > T \gg T_b$, so we can compute $D_s(-j\tau)$  in  the homogeneous medium with permittivity $\eps_o$. 

 It follows from the calculations 
in appendix \ref{ap:A}  that the Toeplitz part of $\mathbb{G}_s$ comes from either the inner product of the primary snapshots 
\eqref{eq:I3} or of the adjoint snapshots~\eqref{eq:I8}.  The Hankel part of $\mathbb{G}_s$ comes from the cross-inner products of the 
primary snapshots with the adjoint snapshots. Since the direction of propagation of electromagnetic waves is along the vector 
product of the electric and magnetic fields, we see from definitions~\eqref{eq:I5} and~\eqref{eq:I8} that the primary and adjoint waves 
propagate in opposite directions. Thus, it is the Hankel part of $\mathbb{G}_s$ that accounts for the products of forward and backward going waves 
i.e., for reflections. It is represented in the data at the sum of travel times of the two waves 
\[
[T + (m+j) \tau] + (T + m \tau) = 2 T + (2m+j) \tau.
\]  
The Toeplitz part of $\mathbb{G}_s$ accounts for waves propagating in the same direction, starting from the reference time $T$.
It is represented in the data at the difference of travel times 
\[
[T + (m+j) \tau] - (T + m \tau) = j \tau.
\]  

Typically,   $M \tau \ll T$ in SAR imaging. Recalling~\eqref{eq:I2}, this means that   the time of travel of the waves in the heterogeneous and unknown part of the medium  is small with respect to the travel time $T$ in the homogeneous medium, between the antenna and the imaging region. Therefore, we deduce from equation~\eqref{eq:I11} that   the information about the unknown $\eps(\bx)$ is contained entirely in the Hankel part of $\mathbb{G}_s$. 

\vspace{-0.1in}
\subsection{Factorization of the snapshots}
 \label{sect:factU}

Gather the snapshots~\eqref{eq:snapU} in the row vector field 
\begin{equation}
\bU_s(\bx) = \left(u_{s,0}(\bx), \ldots, u_{s,M-1}(\bx) \right),
\label{eq:I12}
\end{equation}
and denote by $\mbox{range}[\bU_s(\bx)]$  the $M$-dimensional space spanned by them. 
This space is not known  in the inverse problem, because we cannot measure inside the medium and we cannot compute 
$\bU_s(\bx)$ without knowing $\eps(\bx)$.  However, it turns out that there is key information contained in $\bU_s(\bx)$ that can be computed
from the data, as we now explain:

Let us factorize  $\bU_s(\bx)$  using the Gram-Schmidt orthogonalization of its components \cite{Golub}
\begin{equation}
\bU_s(\bx) = \bV_s(\bx) \bR_s.
\label{eq:I14}
\end{equation}
This gives the orthonormal basis of  $\mbox{range}[\bU_s(\bx)]$, stored in the $M$-dimensional row vector field 
\begin{equation}
\bV_s(\bx) = \left(v_{s,0}(\bx), \ldots, v_{s,M-1}(\bx) \right), 
\label{eq:I13}
\end{equation}
satisfying 
\begin{equation}
 \int_{\RR^2}\hspace{-0.05in} d \bx \, \bV_s^T(\bx) \bV_s(\bx) = \bI_{M},
 \label{eq:I14Ort}
 \end{equation}
where $\bI_{M} \in \RR^{M \times M}$ is the identity matrix. 
The second factor in~\eqref{eq:I14} is an upper  triangular matrix  $\bR_s \in \RR^{M \times M}$.  

The factorization~\eqref{eq:I14}  is causal, because   
\begin{equation}
v_{s,m}(\bx) \in \mbox{span} \{ u_{s,0}(\bx), \ldots, u_{s,m}(\bx)\}, 
\end{equation}
and conversely,
\begin{equation}
u_{s,m}(\bx) \in \mbox{span} \{ v_{s,0}(\bx), \ldots, v_{s,m}(\bx)\},
\end{equation}
for all $0 \le m \le M-1$. 
As was the case with $\bU_s(\bx)$, we cannot compute $\bV_s(\bx)$ without knowing $\eps(\bx)$. 
However,  the upper triangular  $\bR_s$ 
can be computed from the data, because it is, in fact, the Cholesky square root of $\mathbb{G}_s$ that is known by Theorem~\ref{prop.2}:
\begin{align}
\mathbb{G}_s &= \int_{\RR^2} d \bx \, \bU_s(\bx)^T \bU_s(\bx) \nonumber \\
&\stackrel{\eqref{eq:I14}}{ =} \bR_s^T \int_{\RR^2} d \bx \, \bV_s^T(\bx) \bV_s(\bx) \bR_s \stackrel{\eqref{eq:I14Ort}}{ =} \bR_s^T \bR_s. \label{eq:I17}
\end{align}

\vspace{-0.15in}
\subsection{Approximation of the snapshots}
 \label{sect:approxU}

We deduce from Theorem \ref{prop.2} that we have a linear and bijective mapping between $\mathbb{G}_s$ and the data set
\begin{equation}
\label{eq:I16}
%\hspace{-0.01in}
\Big\{D_s(t), t \in \{j \tau, 2T +  (j+m) \tau\}, 0 \le j, m \le M-1 \Big\}.
\end{equation}
The existence and uniqueness of the  Cholesky factorization~\eqref{eq:I17} implies that 
the map between the data~\eqref{eq:I16} and the block upper triangular matrix $\bR_s$  is also bijective. It is $\bR_s$ that ensures the data fit. The uncomputable basis $\bV_s(\bx)$ plays  no role in the expression of $\mathbb{G}_s$ and therefore in the data fit.

Motivated by this observation, we introduce  a causal family of ``internal wave" snapshots, parametrized by the search permittivity $\hat{\eps}(\bx)$. These are the $M$ components of the  row vector field\footnote{Note our notation convention: When the operator, the fields and matrices correspond to the search  permittivity $\hat{\eps}(\bx)$, we indicate it in the arguments. If the permittivity is the true and unknown one, we drop the  argument.} 
\begin{equation}
\tilde{\bU}_{s}(\bx;\hat{\eps}) = \bV_s(\bx;\hat{\eps}) \bR_s, 
\label{eq:I18}
\end{equation}
defined by the orthonormal basis  $\bV_s(\bx;\hat{\eps})$  of the space $\mbox{range}[\bU_s(\bx;\hat{\eps})]$. This satisfies the Gram-Schmidt 
equation \begin{equation}
\bU_s(\bx;\hat{\eps}) = \bV_s(\bx;\hat{\eps}) \bR_s(\hat{\eps}),
\label{eq:I19}
\end{equation}
where the left hand side is obtained from the analogues  of equations~\eqref{eq:I3},~\eqref{eq:I4_s} and~\eqref{eq:snapU},  with $\cL$ replaced by $\cL(\hat{\eps})$, defined by $\hat{c}(\bx) = 1/\sqrt{\mu \hat{\eps}(\bx)}$ instead of $c(\bx)$.

The difference between equations~\eqref{eq:I18} and~\eqref{eq:I19} is that 
$\bR_s(\hat{\eps})$ is the Cholesky square root of the Gramian computed from the synthetic data predicted by $\cF(\hat{\eps})$. These are not the true data so $\bR_s(\hat{\eps}) \ne \bR_s$ and consequently, $ \tilde{\bU}_{s}(\bx;\hat{\eps}) \ne \bU_s(\bx;\hat{\eps})$.

The internal waves are causal, because if we take the $m^{\rm th}$ column in equation~\eqref{eq:I18}, we have 
\begin{align}
\tilde{u}_{s,m}(\bx) &\in \mbox{span}\{ v_{s,0}(\bx;\hat{\eps}), \ldots , v_{s,m}(\bx;\hat{\eps})\} \nonumber \\
&\stackrel{\eqref{eq:I19}}{=} 
\mbox{span}\{ u_{s,0}(\bx;\hat{\eps}), \ldots , u_{s,m}(\bx;\hat{\eps})\},
\end{align}
for all $0 \le m \le M-1$. They also fit the data by construction, because as we had in~\eqref{eq:I17}, this time with the 
orthonormal basis $\bV_s(\bx;\hat{\eps})$, we get
\begin{align}
\int_{\RR^2} \hspace{-0.05in} d \bx \, \tilde{\bU}^T_{s}(\bx;\hat{\eps})\tilde{\bU}_{s}(\bx;\hat{\eps}) &= \bR_s^T \bR_s = \mathbb{G}_s.
\label{eq:I20}
\end{align}
However, the components of~\eqref{eq:I18} are not snapshots of solutions of the wave equation, unless $\hat{\eps}(\bx) = \eps(\bx)$. 
To drive the search permittivity to the true one, we formulate in the next section an iterative optimization procedure that penalizes the misfit 
between the internal waves and the snapshots computed at  $\hat{\eps}(\bx)$, while keeping, by default,  the data fit.
%
%Our inversion method combines the optimization over all the slow times and minimizes over $\hat{\eps}(\bx)$ the objective function
%\begin{equation}
%\mathcal{O}(\hat{\eps}) = \sum_{s=1}^{S} \| \bI_{M+1} - \bR_s(\hat{\eps}) \bR_s^{-1}\|_F^2,
%\label{eq:I22}
%\end{equation}
%where the terms in the sum are a modification of the right hand side of~\eqref{eq:I21}, with the same global minimum. In the next section we 
%assess the performance of the optimization~\eqref{eq:I22} with numerical simulations and we compare 
%it with that of the data fit optimization~\eqref{eq:2}.

\section{Computational setup}
\label{sect:numerics}
First, we specify the excitation. 
Then, we introduce a modified approximation of the snapshots, that is equivalent to that in section~\ref{sect:approxU}, but
allows faster computations. We end with the formulation of the inverse problem as an optimzation.

\vspace{-0.15in} 
\subsection{The probing beam}
Here we use the rotated system of coordinates $\bx_s = (x_s,y_s)$.
The  forcing is of the form~\eqref{eq:S2}, with  
\begin{equation}
b(x_s,t) = \frac{c_o}{2 \pi} \int_{\RR} dk \, \Big[ e^{-i k c_o t- \frac{x_s^2}{2 r_0^2} +  \frac{i q_0 x_s^2}{2}} f[c_o (k-k_o)] + \mbox{c.c.} \Big]
\label{eq:defbeam}
\end{equation}
and Gaussian 
%$ f(\om) = \exp[-\om^2 (3 T_b)^2/2]$, 
{$ f(\om) = \exp[-\om^2 (T_b/3)^2/2]$}, 
with variance chosen so that 
%$ b(x_s,t) \approx 0$ if  $|t| > T_b$. 
{$ |b(x_s,t)| / \max_{t'} |b(x_s,t')| \leq \exp(-9/2) \approx 0$ if  $|t| > T_b$. }
Here $k = \om/c_o$ is the wave number at frequency $\om$, $c_o = 1/\sqrt{\mu \eps_o}$ is the reference wave speed and $k_o = \om_o/c_o$. The initial radius 
$r_0$ of the beam is commensurate with the size of the phased array modeling the antenna, and 
$q_0$ is the initial quadratic phase.   The ``$+$c.c." denotes the addition of the complex conjugate of the first term, so that 
the forcing is real. 

Our choice of the source profile gives a Gaussian beam 
in the paraxial scaling regime, where at range scale $y_s$, we have
\begin{equation}
k_o \gg \big[ r_0^{-4} + q_0^2\big]^{1/4} \gg 1/y_s.
\label{eq:N3b}
\end{equation}
%If we use the rotated system of coordinates $\bx_s = (x_s,y_s)$, t
The first component of the wave at $t < T$ is   \cite[Chapter 4]{AndrewsPhillips}
\begin{align}
&u_s(\bx_s,t) = \int_{\RR} \frac{d k}{4 \pi} \Big\{{f} [c_o(k-k_o)]  \Big(\frac{r_0}{r_{y_s}}\Big)^{1/2} e^{-i k c_o t} \nonumber \\
& \hspace{0.5in} \times e^{ -\frac{x_s^2}{2 r_{y_s}^2} +  \frac{i k x_s^2}{2 \chi_{y_s}} - \frac{i}{2} \arctan \Big(\frac{y_s}{L_R}\Big) + i k y_s} + \mbox{c.c.} \Big\}.
\label{eq:N4}
\end{align}
Here  $L_R = kr_0^2$ is the Rayleigh length,  while the beam radius $r_{y_s}$  and the curvature radius $\chi_{y_s}$ are defined by
\begin{align}
r_{y_s} &= r_0 \sqrt{ \Big(1 + \frac{q_o y_s}{k} \Big)^2 + \frac{y_s^2}{L_R^2}},
\label{eq:N5} \\
\chi_{y_s} &= {(r_{y_s}/r_0)^2}\Big[{\frac{y_s}{L_R^2} + \frac{q_0}{k} \Big(1 + \frac{q_o y_s}{k} \Big)}\Big]^{-1} .
\label{eq:N6}
\end{align}
Note that since 
$
\frac{r_{y_s}}{y_s} \stackrel{y_s \to \infty}{\longrightarrow} r_0 \sqrt{\frac{q_0^2}{k^2} + \frac{1}{L_R^2}},
$
the beam radius without the quadratic phase $q_0$ is smaller than the beam with a quadratic phase at very long range. 
However, $q_0$ is useful for focusing the beam at moderate range. Indeed, equation~\eqref{eq:N5} shows that if $q_0 < 0$, 
the radius $r_{y_s}$ decays with $y_s$ up to the range (focal length)
$
L_\star = {|q_0|k} [{q_0^2 + r_0^{-4}} ]^{-1}
$,
where the beam has minimal radius (beam waist) 
$
r_\star = r_0 [ 1 + q_0^2 r_0^4]^{-1/2}
$.
To have a well focused beam that probes the imaging region at the range of order $y_s$, 
we can choose the quadratic phase $q_0 \approx -k/y_s$ so that $L_\star \approx y_s$
and  $r_\star \approx y_s/(k r_0)$.

\vspace{-0.1in}
\subsection{The internal waves}
The inversion procedure decribed in the next section is iterative and requires repeated evaluations 
of the Gramian $\mathbb{G}_s(\hat{\eps})$ and its square root $\bR_s(\hat{\eps})$, for $\hat{\eps}(\bx)$ updated at each iteration. 
These involve just the first component  of the  wave.  Thus, instead of dealing with the first-order system,   it is computationally advantageous 
to solve the second-order scalar wave equation 
\begin{align}
\partial_t^2 u_s(\bx,t;\hat{\eps}) +A(\hat{\eps})u_s(\bx,t;\hat{\eps})  = 0, \quad t > T  \label{eq:N3}
\end{align}
with initial conditions 
\begin{align*}
&u_s(\bx,T;\hat{\eps}) = u_{s,0}(\bx) = \mbox{ eq.~\eqref{eq:N4} evaluated at } t = T \\
&\partial_t u_s(\bx,T;\hat{\eps}) = \mbox{ time derivative of~\eqref{eq:N4} evaluated at } t = T.
\end{align*}
Equation~\eqref{eq:N3} is deduced from ~\eqref{eq:S4}, with the wave speed $1/\sqrt{\mu \hat{\eps}(\bx)}$ and with the self-adjoint and positive definite 
operator
\begin{equation}
A(\hat{\eps}) = - \frac{1}{\sqrt{\mu \hat \eps(\bx)}} \Delta \Big[\frac{1}{\sqrt{\mu \hat \eps(\bx)}} \cdot \Big].
\label{eq:defA}
\end{equation} 
When $\hat{\eps}(\bx)$ equals the true $\eps(\bx)$, the solution  evaluated at the  instances \eqref{eq:I1} gives the snapshots that we wish to approximate. But no matter what $\hat{\eps}(\bx)$ is, since the support of $\hat{\eps}(\bx)-\eps_o$ is not reached by the waves until $t = T$, we have the same initial conditions at $t = T$ as for the true wave.

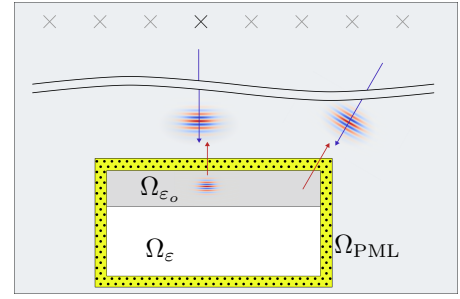
\begin{figure}[t]
\begin{minipage}{0.5\textwidth}
\begin{center}
\begin{tikzpicture}
   % \draw (0, 0) node[inner sep=0] {\includegraphics[width=0.85\textwidth]{Inkskape/drawing.pdf}};
   \draw (0, 0) node[inner sep=0] {\includegraphics[width=0.65\textwidth]{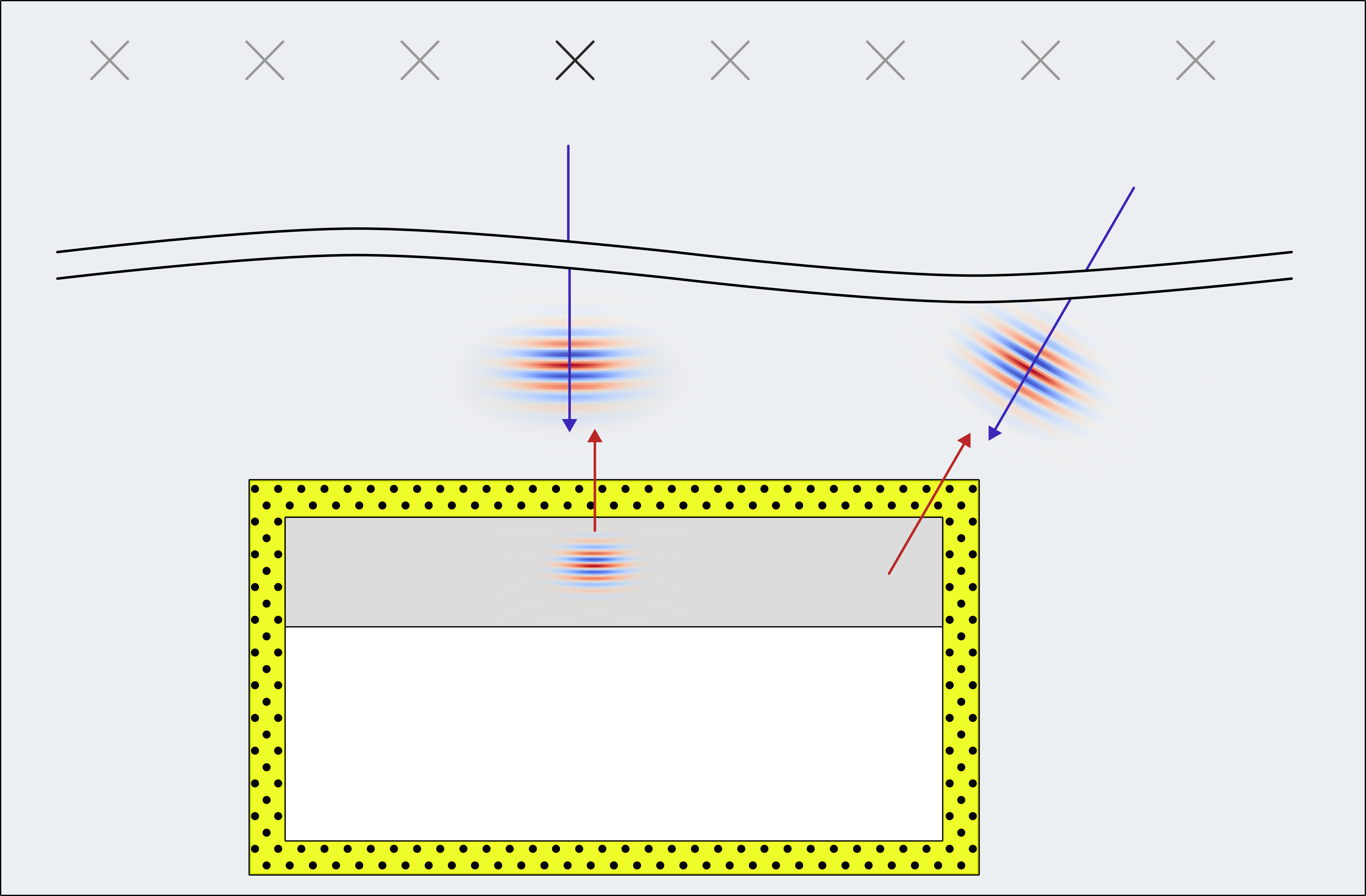}};
\draw (-1., -0.55) node {$\Omega_{\eps_o}$};
\draw (-1., -1.35) node {$\Omega_{\eps}$};
\draw (1.75, -1.3) node {$\Omega_{\rm PML}$};
%\draw [ black,thick] (1.6, -2.2) -- (1.8, -1.9);
%forreferecening
%\draw (-1, -1) node {(-1,-1)};\draw (-2, -2) node {(-2,-2)};\draw (0, 0) node {(0,0)};	\draw (2, 2) node {(2,2)};    \draw (1, 1) node {(1,1)};
\end{tikzpicture}
\end{center}
\end{minipage}
\vspace{-0.05in}\caption{{
Antennas at different locations on the flight track emit beams that enter the computational domain surrounded by a PML boundary layer drawn in yellow, designed to absorb the outgoing waves.
The domain is divided into the $\Omega_{\eps_o}$ part with $\eps(\bx) = \eps_o$ where $u_{s,0}(\bx)$ is supportend and  $\Omega_\eps$ that contains the support of $\eps(\bx)-\eps_o$. }
}
 \label{fig:setup}
\end{figure}

\vskip 0.05in
\subsubsection{Computational domain}
In SAR the distance between the antenna and the imaging region is large, but much of it is contained in the known and homogeneous 
region with permittivity $\eps_o$.  This can be used to speed up the inversion by restricting the computation of  $u_s(\bx,t;\hat{\eps})$ to a smaller domain $\Omega$,  
surrounded by a perfectly matched layer (PML) that absorbs the outgoing waves, as illustrated in Fig. \ref{fig:setup}. The computational domain  is the union of  $\Omega_{\eps_o}$, where 
the permittivity equals  $\eps_o$, and $\Omega_\eps$, that supports $\eps(\bx)-\eps_o$. The initial condition $u_{s,0}(\bx)$ is 
supported in $\Omega_{\eps_o}$. The imaging domain $\Omega_{\rm im}$ is a subset of $\Omega_\eps$.

Because we compute the waves only in $\Omega$, we cannot directly get the matrix $\mathbb{G}_s(\hat{\eps})$.
What we get are the inner products 
\begin{equation}
\mathcal{C}_j(\hat{\eps}):= \lb u_{s,0},u_s(\cdot,j \tau + T;\hat{\eps})\rb, \quad j \ge 0, \label{eq:Sec0}
\end{equation}
because $u_{s,0}(\bx)$ and therefore $u_{s,0}(\bx) u_s(\bx,j \tau + T;\hat{\eps})$ are supported in $\Omega$.
This gives, according to formula~\eqref{eq:I11} evaluated at $m = 0$, see also appendix \ref{ap:A}, 
\begin{align}
\mathcal{C}_j(\hat{\eps})&= \frac{1}{2} [ D_s(j \tau) + D_s(-j \tau)+ D_s(2 T + j \tau;\hat{\eps})],
\label{eq:Sec1}
\end{align}
where the first two terms are independent of $\hat{\eps}(\bx)$, but the last term depends on it. It is defined as in equation 
\eqref{eq:S6}, with the wave field replaced by $u_s(\bx,t;\hat{\eps})$. 

We deduce from~\eqref{eq:Sec1} that  we can determine the Hankel part of 
$\mathbb{G}_s(\hat{\eps})$ from the computed inner products~\eqref{eq:Sec0}, using 
\begin{align*}
&\frac{1}{2} D_s(2 T + (2m+j) \tau;\hat{\eps}) = \mathcal{C}_{2m+j}(\hat{\eps}) \\
& \qquad - \frac{1}{2} [ D_s((2m+j) \tau) + D_s(-(2m+j) \tau)].
\end{align*}
The entries of the Gramian above the diagonal are then
\begin{align}
&[\mathbb{G}_s(\hat{\eps})]_{m,m+j} = \mathcal{C}_{2m+j}(\hat{\eps})  + \frac{1}{2} [ D_s(j \tau) + D_s(-j \tau)] \nonumber \\
& \qquad - \frac{1}{2} [ D_s((2m+j) \tau) + D_s(-(2m+j) \tau)], \label{eq:Sec2}
\end{align}
for $0 \le m \le M-1$ and $0 \le j \le M-1 - m$. 

\vskip 0.05in 
\subsubsection{Computationally advantageous transformation}
Formula~\eqref{eq:Sec2}, while computable, is not convenient, because it involves the terms $D_s$ 
that are independent of $\hat{\eps}(\bx)$. They can be determined from equation~\eqref{eq:N4}, after integration as in~\eqref{eq:S6}.
However, this adds to the computational burden. To avoid this cost, we use the following trick: Instead of seeking to approximate 
the snapshots~\eqref{eq:snapU}, we approximate 
\begin{equation}
%\hspace{-0.1in} 
u_{s}^{\rm ev}(\bx,m \tau) = \frac{1}{2} \left[ u_s(\bx,T+ m \tau) + u_s(\bx,T- m \tau)\right],
\label{eq:Ev1}
\end{equation}
for $0 \le m \le M-1$. The approximations are equivalent, because in the last term in~\eqref{eq:Ev1} the wave is evaluated before 
time $T$ and is thus independent of $\eps(\bx)- \eps_o$.  

Define $u_s^{\rm ev}(\bx,t;\hat{\eps})$ like in~\eqref{eq:Ev1}, by replacing the first term in 
the right-hand side  with $u_s(\bx,m\tau + T;\hat{\eps})$. Note that this wave is even in time (hence the index ``ev") 
and satisfies 
\begin{align}
&\partial_t^2 u_s^{\rm ev}(\bx,t;\hat{\eps}) +A(\hat{\eps})u_s^{\rm ev}(\bx,t;\hat{\eps})  = 0, \quad t > 0 , \label{eq:Ev0}\\
&u_s^{\rm ev}(\bx,0;\hat{\eps}) = u_{s,0}(\bx), \quad \partial_t u_s^{\rm ev}(\bx,0;\hat{\eps}) = 0. \label{eq:Ev01}
\end{align}
In our simulations we solve~\eqref{eq:Ev0}-\eqref{eq:Ev01} in the 
setup in Fig. \ref{fig:setup} (see appendix \ref{ap:B}). Then, we compute the ``new predicted data" 
\begin{align}
\mathbb{D}^{\rm ev}_s(j \tau;\hat{\eps}) &=  \lb u_{s,0}, u_s^{\rm ev}(\cdot,j \tau;\hat{\eps} )\rb, \quad j \ge 0,
\label{eq:S2_C}
\end{align}
that according to the next theorem, proved in appendix \ref{ap:C},  determine the Gramian of the even wave snapshots. Again, we can compute~\eqref{eq:S2_C}  in our setting because the 
product $u_{s,0}(\bx)u_s^{\rm ev}(\bx,j \tau;\hat{\eps} )$ is supported in $\Omega$.

\vskip 0.1in
\begin{thm}
\label{thm.2}
\emph{
Assume that the probing signal is even in time.
Denote by $\mathbb{G}_{s}^{\rm ev}(\hat{\eps})$ the Gramian of the even snapshots at search permittivity
$\hat{\eps}(\bx)$. Its  entries on and above the diagonal are given by 
\begin{align}
\Big[\mathbb{G}_{s}^{\rm ev}(\hat{\eps})\Big]_{m,m+j} &=\lb u_s^{\rm ev}(\cdot,m \tau;\hat{\eps} ),u_s^{\rm ev}(\cdot,(m+j) \tau;\hat{\eps} )\rb \nonumber \\
&\hspace{-0.3in} = \frac{1}{2}  \big[ \mathbb{D}_s^{\rm ev} \big((2m+j)\tau;\hat{\eps}\big) + \mathbb{D}_s^{\rm ev} (j \tau;\hat{\eps})\big],
\label{eq:Ev4}
\end{align}
for $0 \le m \le M-1$ and $0 \le j \le M-1 - m $. The entries below the diagonal are determined by symmetry.
Moreover,  the Gramian  is related to  the one in~\eqref{eq:Sec2} by 
\begin{align} 
[\mathbb{G}^{\rm ev}(\hat{\eps}) ]_{m,m+j} &= \frac{1}{4} \Big([\mathbb{G}_s(\hat{\eps})]_{m,m+j}  + [\mathbb{G}_s(\hat{\eps})]_{0,j} \Big) \nonumber \\
& - \frac{1}{4} \Big([\Lambda_s]_{m,m+j} +  [\Lambda_s]_{0,j}\Big), \label{eq:Ev5}
\end{align} 
where $\Lambda$ is an $M \times M$ symmetric matrix with Toeplitz + Hankel structure,  that is independent of $\hat{\eps}(\bx)$ or $\eps(\bx)$. Its entries on and above the diagonal are 
\begin{align}
\hspace{-0.1in}[\Lambda_s]_{m,m+j} &= 
\frac{1}{2}[D_s(j \tau) + D_s(-j \tau)- D_s(2T - (2m+j) \tau)]\nonumber \\
&-\big[D_s\big((2m+j) \tau\big) + D_s\big(-(2m+j) \tau\big)\big],\label{eq:Ev6}
\end{align}
for $0 \le m \le M-1$ and $0 \le j \le M-1-m$.}
\end{thm}

\vskip 0.1in Since  $\mathbb{G}_s^{\rm ev}(\hat{\eps})$ is completely determined by the computed~\eqref{eq:S2_C},  without any need to subtract 
terms, it is more convenient to work with the even waves~\eqref{eq:Ev1}. The theorem relates 
their Gramian to the one of the snapshots without the even time extension. Since in the end, the approximation of 
\eqref{eq:Ev1} is equivalent to the approximation of the snapshots~\eqref{eq:snapU}, we do not lose any information.

\vskip 0.05in 
\subsubsection{Computation of the internal waves}

The approximation of the even internal waves is analogous to that in section \ref{sect:approxU}: 
\begin{align}
\tilde{\bU}_s^{\rm ev}(\bx;\hat{\eps}) = \bV_s^{\rm ev}(\bx;\hat{\eps}) \bR_s^{\rm ev}.
\label{eq:Ev7}
\end{align}
Here $\bR_s^{\rm ev}$ is the Cholesky square root of the Gramian computed from the measurements, corresponding to 
the true and unknown $\eps(\bx)$,
\begin{align}
\mathbb{G}_s^{\rm ev} = \big(\bR_s^{\rm ev}\big)^T \bR_s^{\rm ev}.
\label{eq:Ev8}
\end{align}
It follows from the calculations in appendix \ref{ap:C} that the analogue of~\eqref{eq:S2_C}, for $\hat{\eps}(\bx)$ replaced by the  unknown, true 
$\eps(\bx)$, can be deduced from the data~\eqref{eq:S6},
\begin{align}
\mathbb{D}_s^{\rm ev}(j \tau) &= \frac{1}{2} \left[ D_s(j \tau) +  D_s(-j \tau)\right]  \nonumber \\
&+ \frac{1}{4} \left[D_s(2T + j \tau) +  D_s(2T - j \tau)\right]. \label{eq:DEv}
\end{align} 
Using this in equation~\eqref{eq:Ev4} we get the data driven expression of $\mathbb{G}_s^{\rm ev}$ and therefore, we can 
compute $\bR_s^{\rm ev}$.

The orthonormal basis stored in $\bV_s^{\rm ev}(\bx;\hat{\eps})$ is defined by the Gram-Schmidt orthogonalization 
of the even snapshots, computed at the search permittivity $\hat{\eps}(\bx)$ and  stored in 
\begin{align}
\bU_s^{\rm ev}(\bx;\hat{\eps})&= \left(u_s^{\rm ev}(\bx,0;\hat{\eps}), \ldots, u_s^{\rm ev}(\bx,(M-1)\tau;\hat{\eps})\right) \nonumber \\
&=\bV_s^{\rm ev}(\bx;\hat{\eps}) \bR_s^{\rm ev}(\hat{\eps}).\label{eq:Ev9}
\end{align}
Here $\bR_s^{\rm ev}(\hat{\eps})$ is the Cholesky square root of the Gramian 
of the even snapshots computed at the search permittivity,
\begin{equation}
\mathbb{G}_s^{\rm ev}(\hat{\eps}) = \big(\bR_s^{\rm ev}(\hat{\eps})\big)^T \bR_s^{\rm ev}(\hat{\eps}).
\label{eq:Ev10}
\end{equation}

\vspace{-0.1in}
\subsection{Optimization formulation of inversion}

Similar to what we explained in section \ref{sect:approxU}, all the information in the data is contained in the 
Cholesky square root $\bR_s^{\rm ev}$ of the Gramian in~\eqref{eq:Ev8}. The orthonormal basis is irrelevant for the data fit. Its purpose in~\eqref{eq:Ev7} 
is to map the upper triangular matrix $\bR_s^{\rm ev}$ to the space of $\bx$-dependent waves that are all consistent with the data. However, these waves are not solutions of the equation~\eqref{eq:Ev0}, unless $\hat{\eps}(\bx) = \eps(\bx)$. This motivates minimizing over $\hat{\eps}(\bx)$ the solution misfit 
\begin{equation}
%\hspace{-0.08in}
\|\tilde \bU_s^{\rm ev}(\cdot;\hat{\eps})- \bU_s^{\rm ev}(\cdot;\hat{\eps})\|^2 = \|\bV_s^{\rm ev}(\cdot;\hat{\eps}) \big[ \bR_s^{\rm ev} - \bR_s^{\rm ev}(\hat{\eps})\big] \|^2  ,\label{eq:Ev11}
\end{equation}
where  
$\| \bU_s(\cdot) \|^2=\sum_{m=0}^{M-1} \int_{\RR^2} |u_{s,m}(\bx)|^2 d\bx$.

%where $\| \cdot \|^2$  is the sum of the squared norms of the components. These norms are  induced by the inner product $\lb \cdot, \cdot \rb$.

Note that in~\eqref{eq:Ev11} we minimize the discrepancy of two row vector fields lying in the same space, with orthonormal basis stored in $\bV_s^{\rm ev}(\bx;\hat{\eps})$.
The Euclidian norm is independent of the basis, so the objective function is in fact the discrepancy of the Cholesky square roots, measured in the squared Frobenius norm: $~\|\bR_s^{\rm ev} - \bR_s^{\rm ev}(\hat{\eps})\|_F^2$, where  
$\| \bR \|_F^2=\sum_{m,j=0}^{M-1}   |R_{m,j} |^2$. 

Our inversion procedure minimizes the objective function
\begin{equation}
\mathcal{O}(\hat \eps) = \sum_{s = 1}^S \|\bI_M - \bR_s^{\rm ev}(\hat{\eps}) (\bR_s^{\rm ev})^{-1} \|_F^2,
\label{eq:Ev12}
\end{equation}
where $S$ is the total number of slow times. The terms in~\eqref{eq:Ev12} differ from what we described above by the multiplication with $(\bR_s^{\rm ev})^{-1}$.  This does not change the global minimum at $\hat \eps(\bx) = \eps(\bx)$,
but it is useful in amplifying the effect of weak echoes in the data. These echoes are reflected in the smaller eigenvalues of  $\mathbb{G}_s^{\rm ev}$  and therefore of $\bR_s^{\rm ev}$.  

In the next section we compare the results given by the minimization of \eqref{eq:Ev12} and of the  data fitting (FWI) objective function 
\begin{equation}
\mathcal{O}^{{\rm FWI}}(\hat \eps) = \sum_{s = 1}^S \sum_{j=0}^{2(M-1)} |\mathbb{D}_s^{\rm ev}(j \tau) -\mathbb{D}_s^{\rm ev}(j \tau;\hat \eps)|^2.
\label{eq:EvFWI}
\end{equation} 
In light of Theorem \ref{thm.2}, data fitting is basically the same as minimizing 
\begin{equation*}
\sum_{s = 1}^S \|\mathbb{G}_s^{\rm ev} - \mathbb{G}_s^{\rm ev}(\hat{\eps})\|_F^2 = \sum_{s = 1}^S \|[\bR_s^{\rm ev}]^T \bR_s^{\rm ev} - 
[\bR_s^{\rm ev}(\hat \eps)]^T\bR_s^{\rm ev}(\hat \eps)\|_F^2.
\end{equation*}

We use a disk shaped  imaging domain $\Omega_{\rm im}$, where $\hat \eps(\bx)$ is parametrized as a linear combination of Gaussian functions centered at the node points $\bz_q$ in a uniform lattice, with equilateral triangular cells of side $h = O(\la_o)$, where $\la_o = 2 \pi/k_o$ is the central wavelength.  
If we let $Q$ be the number of lattice points contained in $\Omega_{\rm im}$, then 
the parametrization is 
\begin{equation}
\hat \eps(\bx) =  \mathcal{E}_{\balpha}(\bx) := \eps_o + \sum_{q = 1}^{Q} \alpha_q e^{ - \frac{\|\bx-\bz_q\|^2}{2 \sigma^2}},
\label{eq:Par1}
\end{equation}
where $\sigma = O(h)$. We denote by $ \mathcal{E}_{\balpha}(\bx)$ the parametrized permittivities, to emphasize 
their dependence on the vector $\balpha = (\alpha_q)_{q=1}^Q$ of coefficients. The goal of the optimization is to determine these coefficients.
We refer to appendix \ref{ap:B} for the details on the minimization of \eqref{eq:Ev12} and \eqref{eq:EvFWI}, including the regularization penalty. The radius of $\Omega_{\rm im}$ and 
the choice of $h$ and $\sigma$ vary with the numerical simulations and are given in the next section.

\vspace{-0.1in}
\subsection{Summary of the inversion method} We can now summarize the steps of our inversion method: 

\vskip 0.06in
\noindent \textbf{Input:} The electric field $E_s(\bx,t)$ at points in the phased array modeling the antenna, 
for $t \in \mathcal{T}$ and $1 \le s \le S$. 

\noindent For all $1 \le s \le S$ do:

\noindent $\bullet$ Compute $D_s(t)$ defined in~\eqref{eq:S6}, using that at points $\bx$ in the antenna, $u_s(\bx,t) = \sqrt{\eps_o/\mu} \, E_s(\bx,t)$.

\noindent $\bullet$ Compute  $\mathbb{D}_s^{\rm ev}(t)$ from $D_s(t)$, using equation~\eqref{eq:DEv}. 

\noindent $\bullet$ Compute the Gramian $\mathbb{G}_{s}^{\rm ev}$, using equation~\eqref{eq:Ev4}. 

\noindent $\bullet$ Compute $\bR_s^{\rm ev}$ using the Cholesky factorization~\eqref{eq:Ev8}. This is the data driven part of the objective function~\eqref{eq:Ev12}.

\noindent $\bullet$ Run the optimization iteration given. At each iterate $\balpha$, 
$\bR_s^{\rm ev}(\mathcal{E}_\balpha)$ is computed as above, from the numerically simulated data defined in~\eqref{eq:S2_C}.

\noindent \textbf{Output:} The estimated permitivity $\mathcal{E}_\balpha(\bx)$, where $\balpha$ is the optimal vector of coefficients 
given by the optimization.

\begin{figure}
\centering 
  \includegraphics[width=0.5\linewidth]{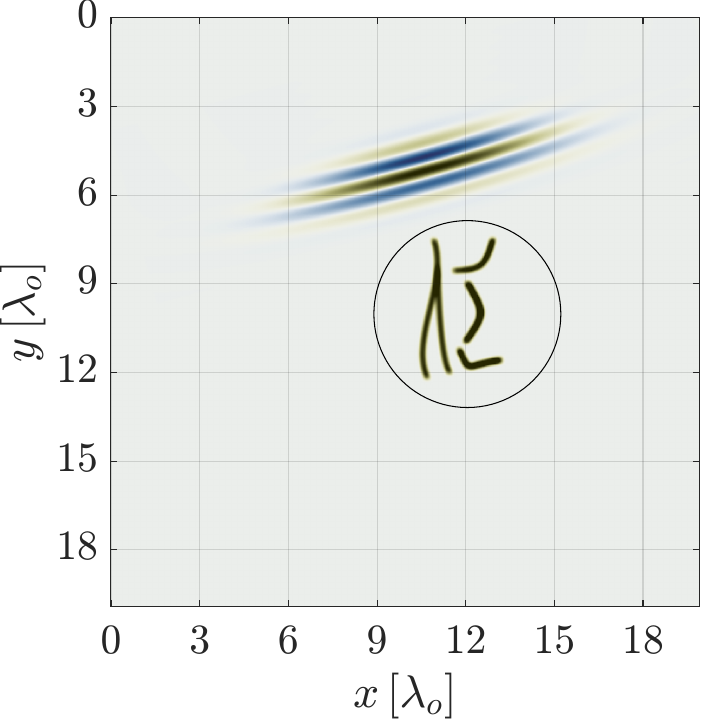}\\
  \caption{Illustration of an initial condition $u_{s,0}(\bx)$ that defines an incident wave beam as it approaches the disk shaped imaging domain $\Omega_{\rm im}$.}
  \label{fig:IW}
\end{figure}

\section{Numerical results}
\label{sect:numerics_res}

The simulations are run in a non-dimensional setting, where the length scales are represented in terms of the central wavelength $\la_o$, 
while the frequencies and bandwidth are in units of $\om_o = k_o c_o$, where $c_o = 3\cdot 10^8$m/s. As an illustration for $X$-band SAR, 
by setting the central frequency to $\om_o/(2 \pi) = 8$GHz, we get the pulse duration $T_b = 0.21$ns and $\la_o = 3.75$cm.
The relative bandwidth at $-3$dB is 66\%.

The initial condition $u_{s,0}(\bx)$ used to solve equations~\eqref{eq:Ev0}-\eqref{eq:Ev01} is computed from equation~\eqref{eq:N4}. 
We use a small, square  computational domain $\Omega$ with side $20 \la_o$. This is because we do iterative optimization, where the cost of computing the objective function over multiple iterations adds up. However, if one is interested in just the first iteration
of the optimization, which as we show below gives a good image of the medium, then it is feasible to enlarge significantly the domain $\Omega$.
The  focused beam 
that defines the initial state $u_{s,0}(\bx)$ has the radius $r_\star = 2.5\la_o$. It corresponds to the excitation \eqref{eq:defbeam} with the antenna of radius 
$r_o = 6.5\,\lambda_o$, 
at range $y_s = 100 \lambda_o$. The quadratic phase is given by $q_0 = -1.1\,\lambda_o^{-2}$.

\vspace{-0.1in} \subsection{First set of simulations}

    \setlength{\figurewidth}{0.44\linewidth}
    \setlength{\figureheight}{0.3\linewidth}
    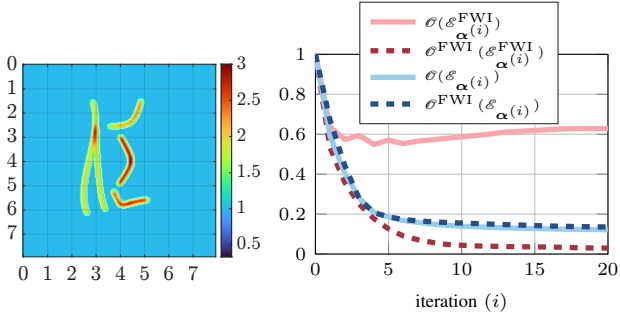
\begin{figure}
        \centering
        \raisebox{15pt}{\includegraphics[width=0.38\linewidth]{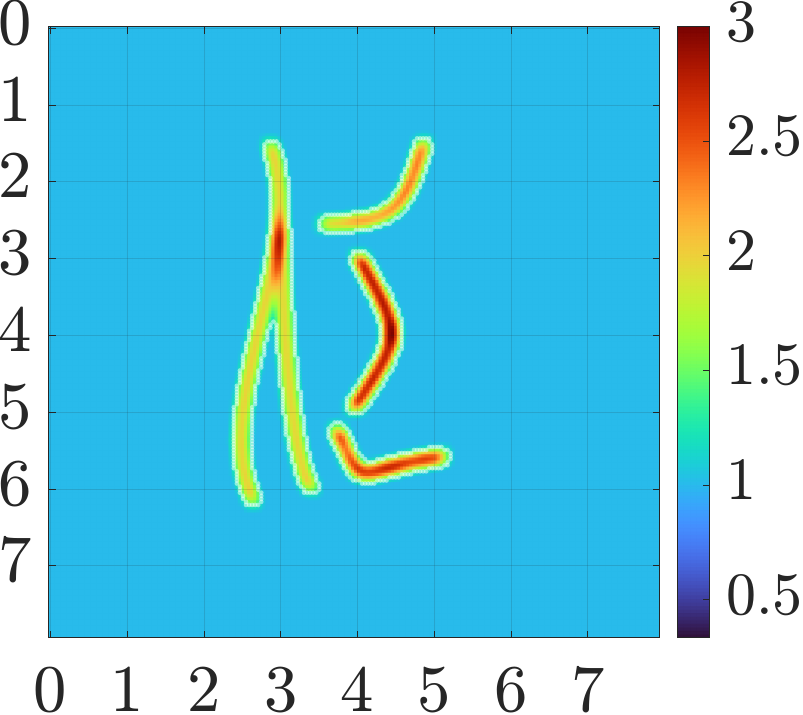}}
\definecolor{mycolor1}{rgb}{0.98000,0.65000,0.68000}%
\definecolor{mycolor2}{rgb}{0.65000,0.20000,0.25000}%
\definecolor{mycolor3}{rgb}{0.60000,0.80000,0.90000}%
\definecolor{mycolor4}{rgb}{0.15000,0.30000,0.55000}%
\begin{tikzpicture}

%\normalsize
%\small
%\footnotesize
%\scriptsize
%\tiny

\begin{axis}[%
width=0.998\figurewidth,
height=\figureheight,
at={(0\figurewidth,0\figureheight)},
scale only axis,
xmin=0,
xmax=20,
xlabel style={font=\tiny, color=white!15!black},
ylabel style={font=\tiny, color=white!15!black},
xlabel={\scriptsize iteration $(i)$},
ymin=0,
ymax=1,
axis background/.style={fill=white},
xmajorgrids,
ymajorgrids,
xlabel near ticks,
tick label style={font=\scriptsize},
label style={font=\tiny},
legend style={at={(0.15,0.65)},font = \tiny, anchor=south west, legend cell align=left, align=left, draw=white!15!black}
]
\addplot [color=mycolor1, line width=2.0pt]
  table[row sep=crcr]{%
0	1\\
1	0.637257779948456\\
2	0.574041793888755\\
3	0.593808936457396\\
4	0.548643181395616\\
5	0.57008005814507\\
6	0.552800427817896\\
7	0.565955732495078\\
8	0.572239610466649\\
9	0.580149533197778\\
10	0.58677889304524\\
11	0.594159955042528\\
12	0.601615252450508\\
13	0.610744755048299\\
14	0.614785245846051\\
15	0.619112871214916\\
16	0.623018902886201\\
17	0.627105628391423\\
18	0.628594846177083\\
19	0.62823335599425\\
20	0.628702191881682\\
};
\addlegendentry{${\cal O}({\cal E}^{\rm FWI}_{\balpha^{(i)}})$}

\addplot [color=mycolor2, dashed, line width=2.0pt]
  table[row sep=crcr]{%
0	1\\
1	0.532175758665597\\
2	0.359750471047821\\
3	0.248919605041492\\
4	0.176929044373958\\
5	0.125791440681612\\
6	0.0909347382329139\\
7	0.0719301019437894\\
8	0.0579699744136357\\
9	0.0481665938024669\\
10	0.0437436847360161\\
11	0.0406867924228276\\
12	0.039147548875197\\
13	0.0373676680365285\\
14	0.0366946494450071\\
15	0.035971035860467\\
16	0.0351593683979374\\
17	0.0336957605894393\\
18	0.0314616193337554\\
19	0.0303376474712794\\
20	0.0299186683220089\\
};
\addlegendentry{$ {\cal O}^{\rm FWI}\, ({\cal E}^{\rm FWI\,}_{\balpha^{(i)}})$}

\addplot [color=mycolor3, line width=2.0pt]
  table[row sep=crcr]{%
0	1\\
1	0.590879530272902\\
2	0.406292370343466\\
3	0.280342830106021\\
4	0.209590779983786\\
5	0.185927003340719\\
6	0.167727336825378\\
7	0.158729569740716\\
8	0.149983575317461\\
9	0.14344174598818\\
10	0.139647099265225\\
11	0.136427139669847\\
12	0.133247804235278\\
13	0.131054617832404\\
14	0.129146054701711\\
15	0.127549907931444\\
16	0.125977027626167\\
17	0.124690652950898\\
18	0.12365051261182\\
19	0.122503313808861\\
20	0.121628142420151\\
};
\addlegendentry{${\cal O}({\cal E}_{\balpha^{(i)}})$}

\addplot [color=mycolor4, dashed, line width=2.0pt]
  table[row sep=crcr]{%
0	1\\
1	0.670911552525686\\
2	0.44432878683629\\
3	0.275303126081378\\
4	0.206542455366177\\
5	0.186361237835923\\
6	0.172675406697861\\
7	0.165975637083907\\
8	0.161437848169658\\
9	0.157807248378587\\
10	0.155369713349778\\
11	0.152747523097154\\
12	0.149949615705927\\
13	0.147392117444818\\
14	0.145219575347324\\
15	0.142851706770286\\
16	0.14044671172139\\
17	0.1387627022231\\
18	0.137502312373895\\
19	0.136539166591829\\
20	0.13543371227104\\
};
\addlegendentry{${\cal O}^{\rm FWI}\,({\cal E}_{\balpha^{(i)}})$}

\end{axis}
\end{tikzpicture}%
        \vspace{-0.15in}\caption{Left:  True permittivity $\eps(\bx)$.  The abscissa and ordinate are in units of $\lambda_o$. The colorbar is in units of $\eps_o$. The probing beams  are incident from the left. Right: Evolutions of the objective functions.}
        \label{fig:True1}
    \end{figure}

In the first simulation we consider thin reflectors, as shown in the left plot of Fig. \ref{fig:True1}. The imaging domain is as shown in Fig. \ref{fig:IW} and 
the triangular mesh size used in equation~\eqref{eq:Par1} to parametrize the search permittivity is {$h = 0.1375 \la_o$}. This gives a vector $\balpha$  of $Q = 2353$ unknown coefficients. We use $S = 31$ slow times, for a circular trajectory of the antenna. The axis of the beam is rotated 
at angles between $-75^\circ$ and $75^\circ$ with respect to the horizontal axis, in equal angle increments of $5^\circ$.
For each slow time, we use $M = 96$ time steps for the computation of the Gramian, at interval $\tau = \frac{\pi}{2.2\om_o}$, corresponding to 2.2 times the Nyquist sampling rate at the central frequency.

The right plot of Fig. \ref{fig:True1} displays convergence curves:  The solid and dashed blue curves show the evolutions of $\mathcal{O}(\mathcal{E}_{\balpha^{(i)}})$ and 
$\mathcal{O}^{\rm FWI}(\mathcal{E}_{\balpha^{(i)}})$, for  the permittivity  $\mathcal{E}_{\balpha^{(i)}}(\bx)$ given by our method. 
The solid and dashed  black curves show the same objective functions computed for the permitivitty 
$\mathcal{E}_{\balpha^{(i)}}^{\rm{FWI}}(\bx)$ given by FWI. Since FWI minimizes the data misfit, it gives a monotone decreasing $\mathcal{O}^{\rm FWI}(\mathcal{E}^{\rm FWI}_{\balpha^{(i)}})$. However, 
 $\mathcal{O}(\mathcal{E}^{\rm FWI}_{\balpha^{(i)}})$ does not decay. Our method gives a slightly worse data fit, but both objective functions improve at each iteration.

\begin{figure}[htbp]
\centering
\begingroup
\setlength{\tabcolsep}{0.3pt}  % horizontal spacing between columns
\renewcommand{\arraystretch}{0}  % vertical spacing between rows

\begin{tabular}{ p{0.32\linewidth} p{0.305\linewidth} p{0.3635\linewidth}}

% --- header row ---
{\centering \textbf{\phantom{FI}Our method}} & {\centering \phantom{FWwI}\textbf{FWI}} & {\centering \phantom{FIi}\textbf{Noisy data}} \\[-0.1mm]

% --- row 1 ---
\includegraphics[width=\linewidth]{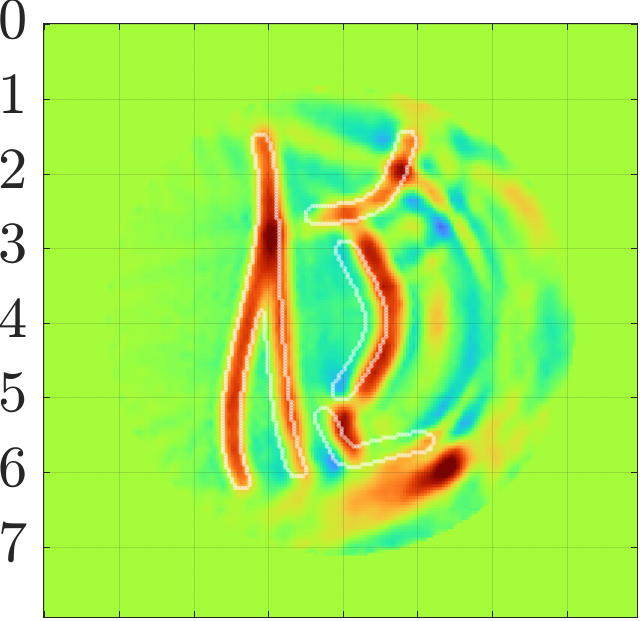} &
\hfill\includegraphics[width=0.98\linewidth]{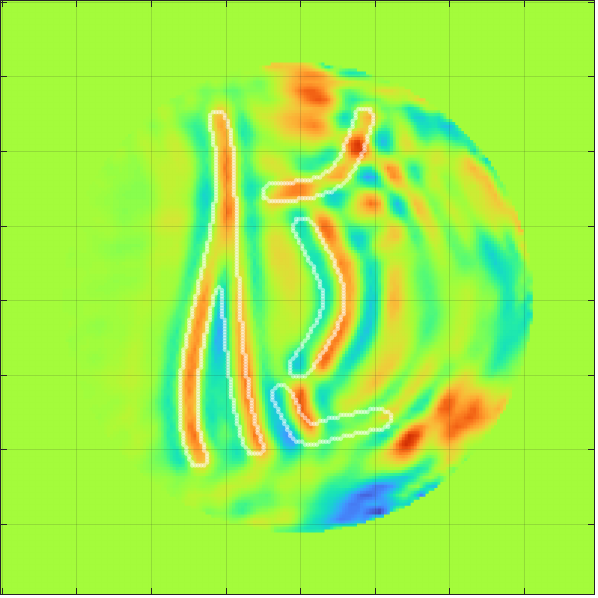} &
\hfill \raisebox{0pt}{\hfill\includegraphics[width=0.9855\linewidth]{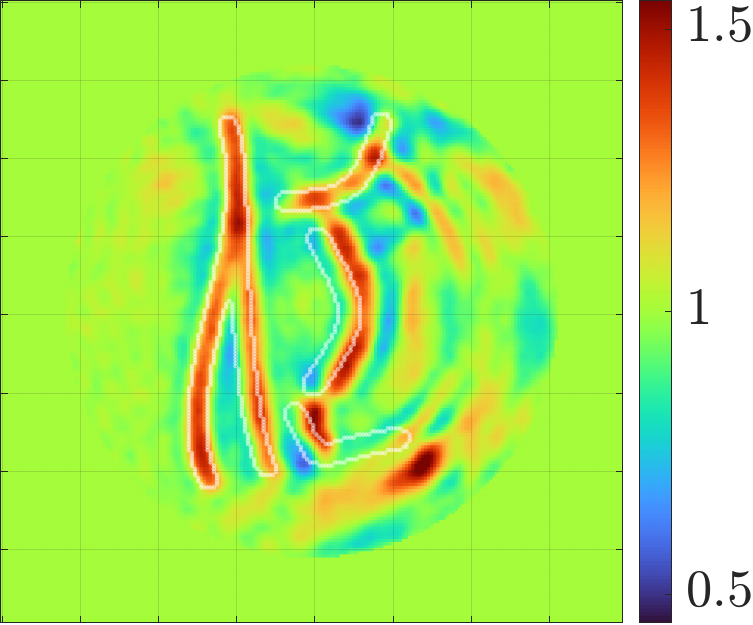}} \\[1mm]

% --- row 2 ---
\includegraphics[width=\linewidth]{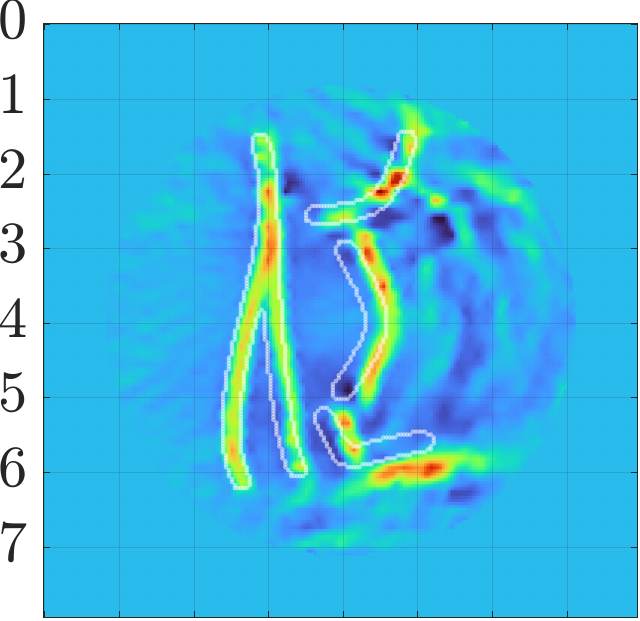} &
\hfill\includegraphics[width=0.98\linewidth]{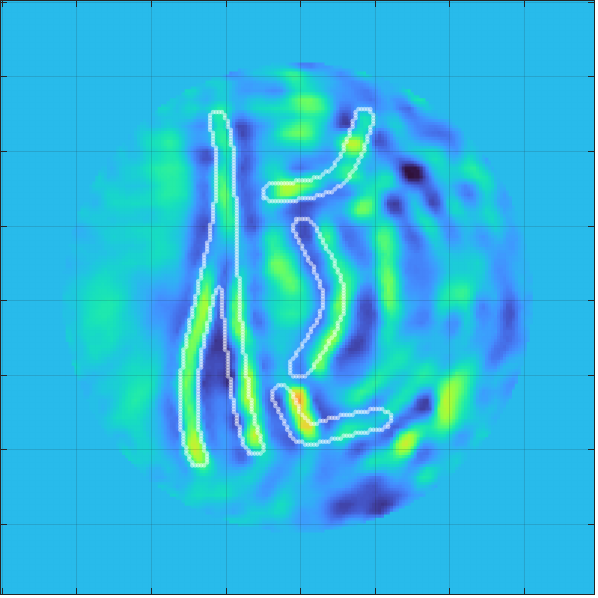} &
\hfill\raisebox{-1.6pt}{\hfill\hspace{1pt}\includegraphics[width=0.9855\linewidth]{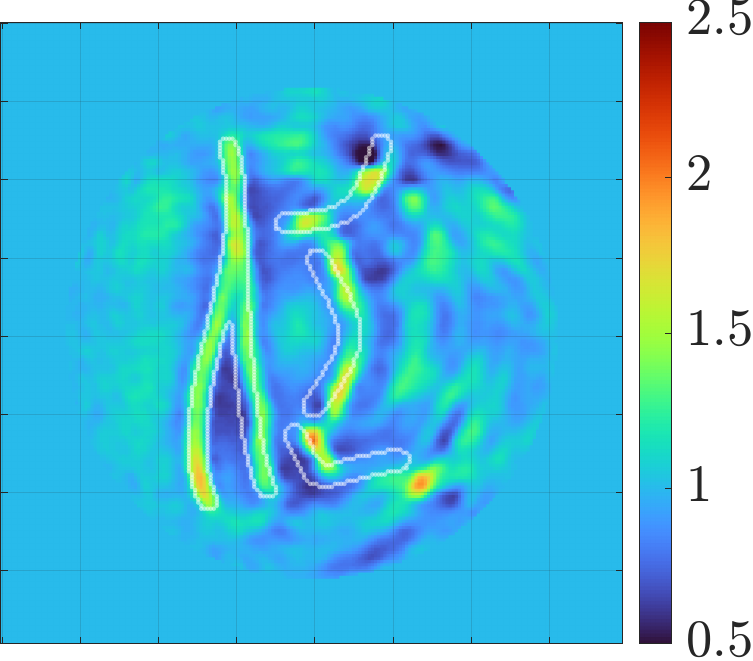}} \\[1mm]

% --- row 3 ---
\raisebox{0.05pt}{\includegraphics[width=\linewidth]{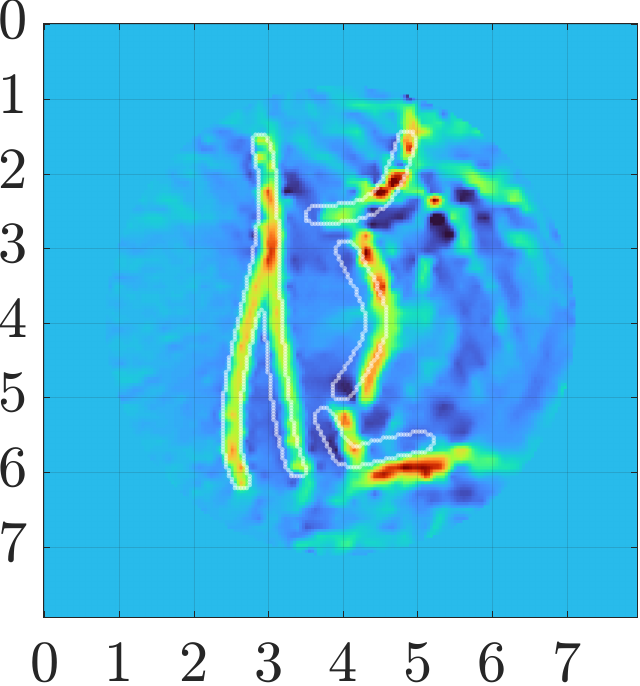}} &
\hfill\includegraphics[width=\linewidth]{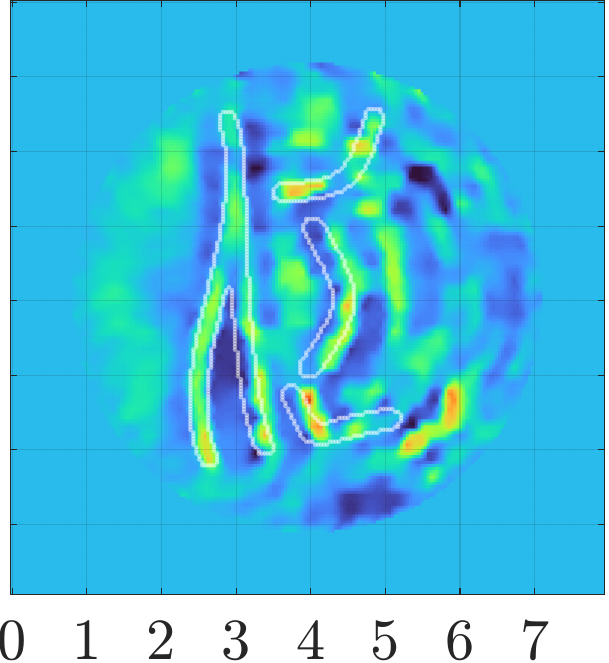} &
\raisebox{.33pt}{\includegraphics[width=\linewidth]{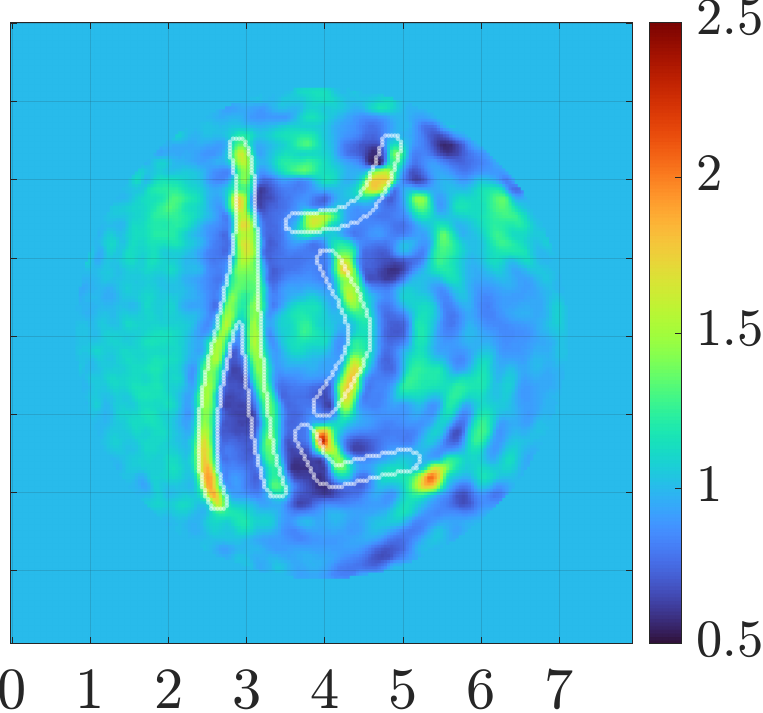}} \\[0.5mm]

\end{tabular}
\endgroup

\caption{{Results of the  inversion: The top row shows the estimated permittivities  after one iteration. The next two  rows show the results after the fifth and tenth iterations, respectively. The left column corresponds to our method, the middle column to FWI, and the right column to our method with noisy data. The white outlines show the true locations of the inclusions. The axes are scaled in multiples of $\lambda_o$. The colorbar is in units of $\eps_o$.}}
\label{fig:Cracks}
\end{figure}

In Fig. \ref{fig:Cracks} we show the inversion results at iterations $i = 1, 5$ and $10$. As seen from the plots in Fig. \ref{fig:True1}, the objective functions change very little after the $10^{\rm th}$ iteration. We superpose the contours of the true inclusions,  to aid in interpreting the results. Note that our method gives a good estimate of the inclusions at the first iteration. It is the contrast that changes at the latter iterations. The FWI approach gives worse estimates of the permittivity at all steps, even though the data are well matched.  The first two columns of plots in Fig. \ref{fig:Cracks} are for noiseless data. The last column shows how our method deals with noise (see appendix \ref{ap:Noise} for the description of the noise model). We do not show FWI results with noisy data, because they do not add any new information and are, naturally, slightly worse than those 
shown in the second column of plots.

\begin{figure}
\centering
\includegraphics[width=0.4\linewidth]{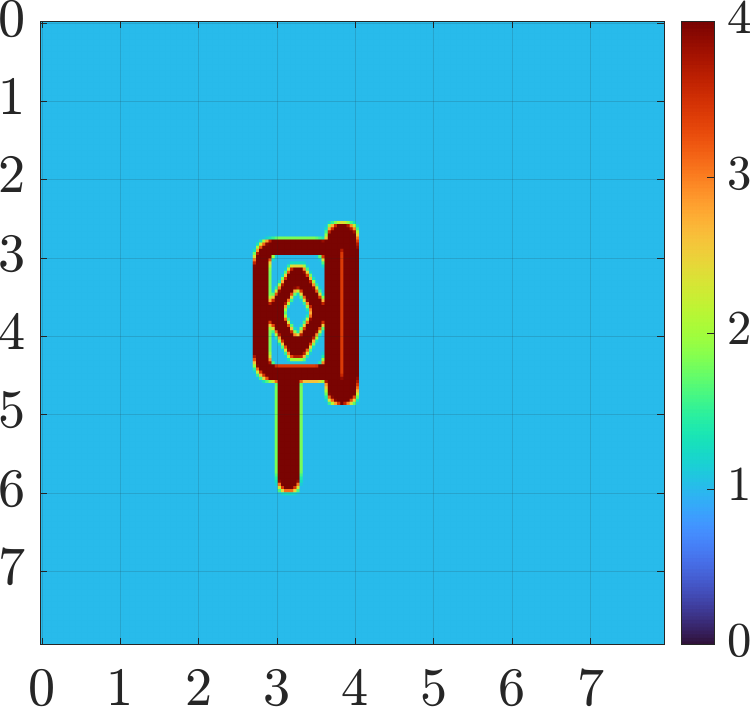}
\caption{True permittivity $\eps(\bx)$ of the second target. The axes are scaled by the center wavelength $\lambda_o$.  The probing beams  are incident from the left. The colorbar is in units of $\eps_o$.}
\label{fig:True2}
\end{figure}

\vspace{-0.1in} \subsection{Second set of simulations}
Here we display results for the more complicated target  shown in Fig. \ref{fig:True2}. The triangular mesh size used in equation~\eqref{eq:Par1} to parametrize the search permittivity is {$h = 0.125 \la_o$}. This gives a vector $\balpha$  of $Q = 2798$ unknown coefficients.  We use $S = 25$ slow times, for a circular trajectory of the antenna. The axis of the beam is rotated 
at angles between $-50^\circ$ and $50^\circ$ with respect to the horizontal axis, in equal angle increments of $4^\circ$.
For each slow time, we use $M = 96$ time steps 
for the computation of the Gramian, at interval {$\tau =  \frac{\pi}{2.2\om_o}$}. 

\begin{figure}
\centering

\begingroup
\setlength{\tabcolsep}{1pt}  % no horizontal space
\renewcommand{\arraystretch}{0}  % no vertical space
\begin{center}
\begin{tabular}{p{0.44\linewidth} p{0.5\linewidth}}

% --- header row ---
{\begin{center}
\textbf{%\phantom{FI}
Our method}\end{center}} &{\begin{center}\textbf{FWI}\phantom{FWI}
\end{center}} \\[-4mm]
% --- row 1 ---
\raisebox{0pt}{\includegraphics[width=0.9\linewidth]{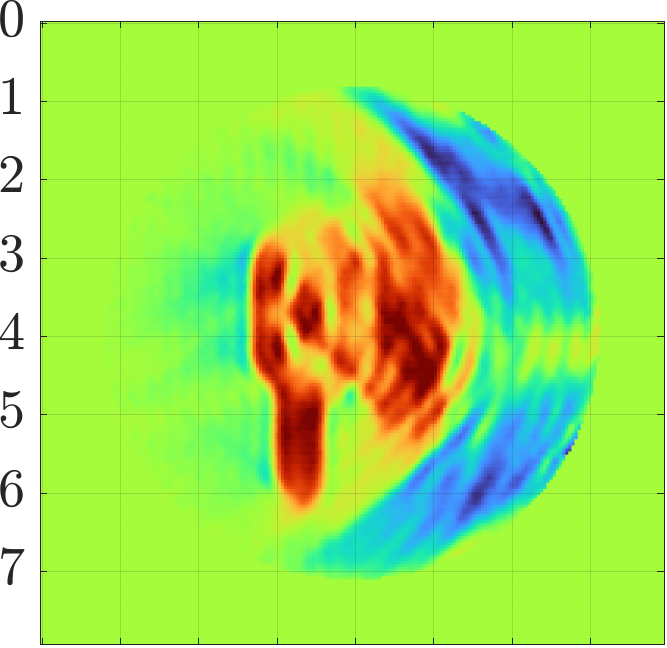}} &
%\hfill
\raisebox{0pt}{\includegraphics[width=0.9\linewidth]{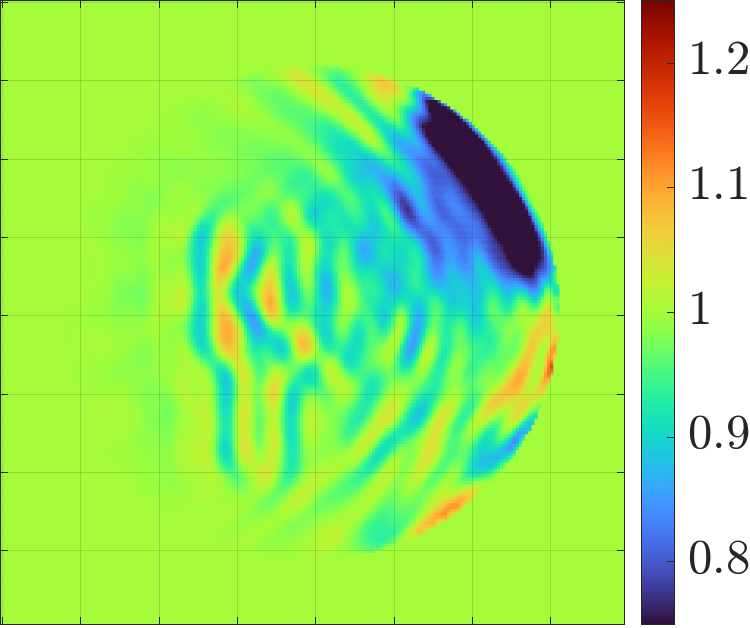}} \\[0.5mm]

% --- row 2 ---
\raisebox{0pt}{\includegraphics[width=0.9\linewidth]{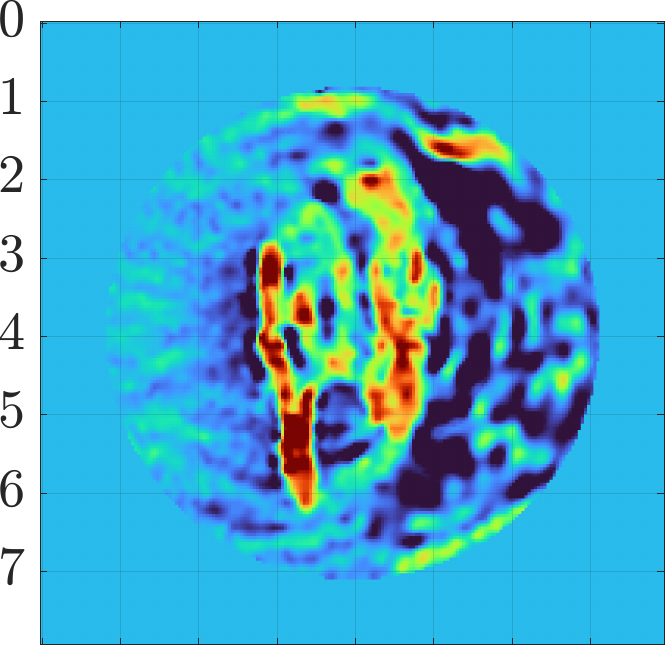}} &
%\hfill
\raisebox{0pt}{\includegraphics[width=0.9\linewidth]{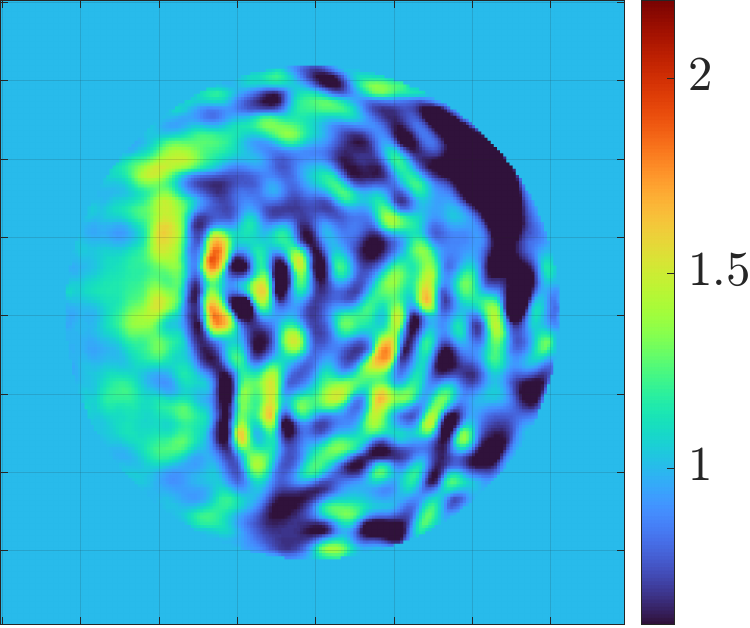}} \\[0.5mm]

% --- row 3 ---
\raisebox{0pt}{\includegraphics[width=0.9\linewidth]{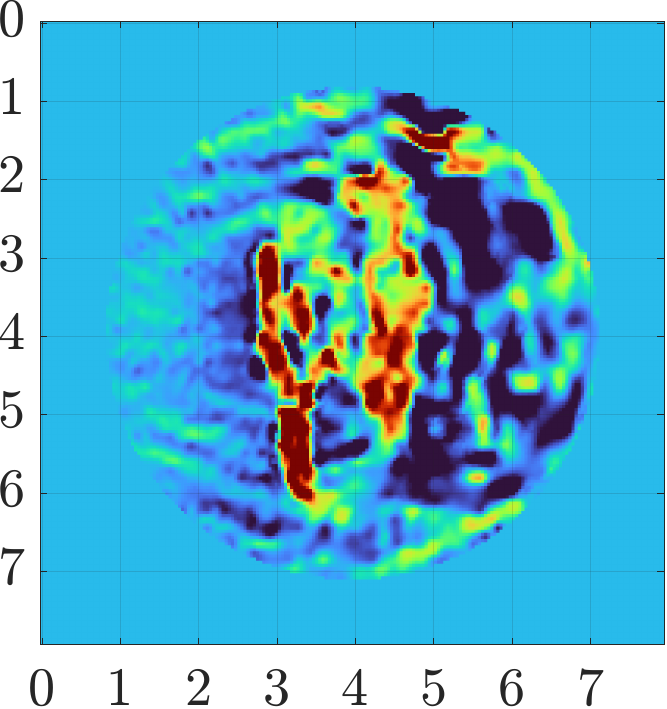}} &
\raisebox{0pt}{\includegraphics[width=0.9\linewidth]{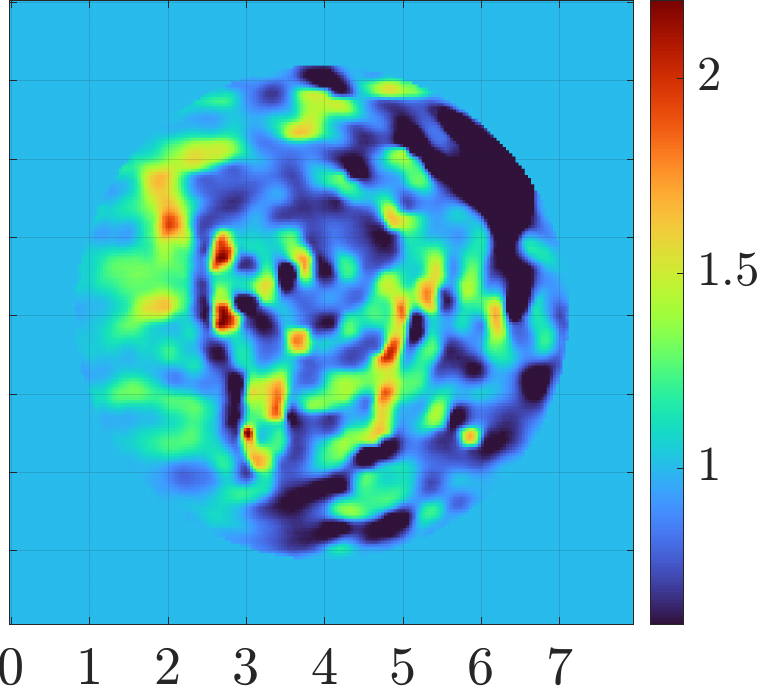}} \\[0.5mm]

\end{tabular}
\end{center}
\endgroup

\caption{Inversion results for the target shown in Fig. \ref{fig:True2}. The first row shows the estimated permitivities after one iteration. The second and third rows show the results after the $10^{\rm th}$ and $30^{\rm th}$ iterations, respectively. The first column corresponds to our method and the right column to FWI. The axes are in units of  $\lambda_o$. The colorbar is in units of $\eps_o$.}
\label{fig:Tanks}
\end{figure}

The convergence curves are similar to those in Fig. \ref{fig:True2}, so we do not include them here. Since the objective functions change very little after the $30^{\rm th}$ iteration, we display in Fig. \ref{fig:Tanks} the inversion results 
at iterations $1, 10$ and $30$. Note again how the first iteration of our method gives a good estimate of the support of the target. The following iterates seek 
to improve the quantitative estimate of the permittivity in the target. All the results are with noiseless data. For brevity, we did not include simulations with noisy data, because they do not bring any new insights.

\section{Summary}
\label{sect:summary}
We introduced a new approach to inverse scattering with synthetic aperture radar (SAR), 
where a moving antenna probes an unknown, heterogeneous, isotropic and nonmagnetic  medium with pulsed, directed beams of radiation and measures the backscattered electric field. The inverse problem is to estimate the dielectric permittivity, and therefore the wave speed in the medium, from the measurements. We showed how to map the measurements to a family of internal waves, at points inside the inaccessible medium. These  waves are computed at the search permittivity. They fit the data by construction, but they do not solve Maxwell's equations unless the permittivity equals the true one. Thus, to estimate the permittivity,  we formulated the inverse problem as an iterative  minimization of the solution misfit,
i.e., the misfit between the internal wave and the solution of Maxwell's equations.
We used numerical simulations to compare the performance of our approach to that of the standard, nonlinear least squares data fitting. We also explained 
that the computational cost of the first iteration of our method is similar to that of standard SAR image formation, but the images are better: they do not have multiple scattering artifacts and give a more accurate estimate of the support of targets.

\section*{Acknowledgment}
This work relates to the  Air Force award number FA9550-22-1-0077, issued by the Air
Force Office of Scientific Research and to the 
Department of Navy award N000142612023, issued by the Office of
Naval Research. It is also partially supported by the Agence de l'Innovation de D\'efense (AID) via
Centre Interdisciplinaire d'\'Etudes pour la D\'efense et la S\'ecurit\'e (CIEDS) project PRODIPO and NSF grant DMS-2309197.
The computations were enabled by resources in project UPPMAX 2025/2-271 provided by Uppsala University at UPPMAX.

\appendices
\section{}
\label{ap:A}
To prove Theorem \ref{prop.2}, we begin with the identity 
\begin{equation}
\int_{\RR^2} \hspace{-0.05in} d \bx \, \bphi_{s,0}^T(\bx) \bphi_{s,j}(\bx) = D_s(j \tau) + D_s(-j \tau).
\label{eq:A1}
\end{equation}
Here we integrate over $\RR^2$, but since we are interested in bounded $j$ and the waves propagate at finite speed, the integrand 
is compactly supported in $\RR^2$.

The derivation of~\eqref{eq:A1} is  in \cite[Section 2]{ROM4} and we summarize it here: Consider 
\begin{equation*}
Q(t;t') = \int_{\RR^2}  \hspace{-0.05in} d \bx \, \left[ \begin{pmatrix} u_s(\bx, t+T) \\ {\bH}_s(\bx,t+T) \end{pmatrix}\right]^T 
\begin{pmatrix} u_s(\bx,t' + t+T) \\ {\bH}_s(\bx,t'+t +T) \end{pmatrix}  
\end{equation*}
and observe from definition~\eqref{eq:I3} that the left-hand side in~\eqref{eq:A1} equals 
$Q(0;j \tau)$. Taking the $t$ derivative and using the system~\eqref{eq:S4} we get 
\begin{align*}
\partial_t Q(t;t') &=  \int_{\RR^2}  \hspace{-0.05in} d \bx \, \left[
\begin{pmatrix}J_s(\bx,t+T) \\ {\bf 0} \end{pmatrix}\right]^T \begin{pmatrix} u_s(\bx,t' + t+T) \\ {\bH}_s(\bx,t'+t +T) \end{pmatrix} \\
&+\int_{\RR^2}  \hspace{-0.05in}d \bx \, \left[ \begin{pmatrix} u_s(\bx, t+T) \\ {\bH}_s(\bx,t+T) \end{pmatrix}\right]^T  \begin{pmatrix}J_s(\bx,t'+t+T) \\ {\bf 0} \end{pmatrix}
\end{align*}
where the two terms involving $\cL$ cancel because $\cL$ is skew-adjoint. Rewriting  the equation above  in terms of the inner product $\lb \cdot, \cdot \rb$ 
\begin{align*}
\partial_t Q(t;t') &=  \lb J_s(\cdot,t+T), u_s(\cdot,t' + t+T)\rb \\
&+\lb u_s(\cdot, t+T),J_s(\cdot,t'+t+T)\rb,
\end{align*}
and integrating from $t = - \infty$ to $t = 0$, we get 
\begin{align*}
Q(0;t') &=  \int_{-\infty}^0 d t  \lb J_s(\cdot,t+T) , u_s(\cdot,t' + t+T)\rb \\
&+ \int_{-\infty}^0 d t \lb u_s(\cdot, t+T), J_s(\cdot,t'+t+T)\rb.
\end{align*}
Change variables in the first term as 
$
t+T \leadsto - t'',
$
and in the second term as 
$
t'+t + T \leadsto - t''.
$
Using that $J_s(\bx,t)$ is supported at $t \in (-T_b,T_b)$ and that $T > T_b$, we get 
\begin{align*}
Q(0;t') &=  \int_{-T_b}^{T_b} dt'' \int_{\RR^2}  \hspace{-0.05in}d \bx \, J_s(\bx,-t'') u_s(\bx,t' -t'') \\
&+ \int_{-T_b}^{T_b} d t" \int_{\RR^2}  \hspace{-0.05in} d \bx \, J_s(\bx,-t'') u_s(\bx,-t' -t'')  \\
&\stackrel{\eqref{eq:S6}}{=} D_s(t')+ D_s(-t'),
\end{align*}
which proves~\eqref{eq:A1}.

Next, we prove that  
\begin{align}
H(j \tau)&:= \int_{\RR^2}  \hspace{-0.05in}d \bx \, \left[\bphi_{s,0}^\star(\bx)\right]^T  \bphi_{s,j}(\bx) \nonumber \\
&= \int_{-T_b}^{T_b} dt' 
\lb J_s(\cdot,t'), u_s(\cdot ,2 T + j \tau - t')\rb \nonumber \\
&=: {\mathcal D}_s(2 T + j \tau).
\label{eq:A2}
\end{align}
Note that $\mathcal{D}_s$ differs from $D_s$, because the convolution is done with the signal without time reversal. If the signal is even in time, as assumed in Theorem \ref{prop.2}, the two are the same. 

Equation~\eqref{eq:A2} follows from 
\begin{align*}
H(t) &= \int_{\RR^2}  \hspace{-0.05in}d \bx \, \left[\bphi_{s,0}^\star(\bx)\right]^T  \bphi_{s}(\bx,t) \\
&\stackrel{\eqref{eq:I4_s}}{=} \int_{-T_b}^{T_b} dt' \int_{\RR^2}  \hspace{-0.05in}d \bx \, \left[ e^{(T-t')\cL} \begin{pmatrix}J_s(\bx,t') \\ {\bf 0} \end{pmatrix}\right]^T
\bphi(\bx,t) \\
&= \int_{-T_b}^{T_b} dt' \int_{\RR^2}  \hspace{-0.05in}d \bx \, \begin{pmatrix}J_s(\bx,t') \\ {\bf 0} \end{pmatrix}^T e^{-(T - t')\cL} \bphi(\bx,t) 
\\
&\stackrel{\eqref{eq:I3}}= \int_{-T_b}^{T_b} dt' \int_{\RR^2}  \hspace{-0.05in} d \bx \, \begin{pmatrix}J_s(\bx,t') \\ {\bf 0} \end{pmatrix}^T \bphi(\bx,t + T-t') \\
&\stackrel{\eqref{eq:I3}}= \int_{-T_b}^{T_b} dt' \lb J_s(\cdot,t'),u_s(\cdot,t + 2 T - t')\rb,
\end{align*}
once we set $t = j \tau$.

The identities~\eqref{eq:A1}-\eqref{eq:A2} are now used to prove the theorem. We have from definitions~\eqref{eq:I9} 
and~\eqref{eq:I10} that 
\begin{align*}
(\mathbb{G}_s)_{m,m+j} &= \frac{1}{4} \int_{\RR^2}  \hspace{-0.05in}d \bx \, [\bphi_{s,m}(\bx)]^T \bphi_{s,m+j}(\bx) \\
 &+ \frac{1}{4} \int_{\RR^2}  \hspace{-0.05in}d \bx \, [\bphi_{s,m}^\star (\bx)]^T \bphi_{s,m+j}^\star (\bx) \\
 &+  \frac{1}{4} \int_{\RR^2}  \hspace{-0.05in}d \bx \, [\bphi_{s,m}^\star (\bx)]^T \bphi_{s,m+j}(\bx) \\
 &+ \frac{1}{4} \int_{\RR^2}  \hspace{-0.05in}d \bx \, [\bphi_{s,m} (\bx)]^T \bphi_{s,m+j}^\star (\bx).
 \end{align*}
 Due to the relation~\eqref{eq:I8} between the primary and adjoint waves, we note that the first two terms are the same and so are the 
 last two terms. Therefore, 
\begin{align*}
(\mathbb{G}_s)_{m,m+j} &= \frac{1}{2} \int_{\RR^2}  \hspace{-0.05in}d \bx \, [\bphi_{s,m}(\bx)]^T \bphi_{s,m+j}(\bx) \\
 &+ \frac{1}{2} \int_{\RR^2}  \hspace{-0.05in}d \bx \, [\bphi_{s,m}^\star (\bx)]^T \bphi_{s,m+j}(\bx).
 \end{align*}
 The first term in this equation is
 \begin{align*}
&\frac{1}{2} \int_{\RR^2}  \hspace{-0.05in}d \bx \, [\bphi_{s,m}(\bx)]^T \bphi_{s,m+j}(\bx) 
\\
&\quad \stackrel{\eqref{eq:I3}}{=}
\frac{1}{2} \int_{\RR^2}  \hspace{-0.05in}d \bx \, \bphi_{s,0}^T(\bx) e^{m \tau \cL} \bphi_{s,m+j}(\bx) \\
&\quad = \frac{1}{2} \int_{\RR^2}  \hspace{-0.05in}d \bx \, \bphi_{s,0}^T(\bx) \bphi_{s,j} (\bx) \\
&\quad \stackrel{\eqref{eq:A1}}{=} \frac{1}{2} \left[ D_s(j \tau) + D_s(-j \tau) \right].
\end{align*}
The second term is 
\begin{align*}
&\frac{1}{2} \int_{\RR^2}  \hspace{-0.05in}d \bx \, [\bphi_{s,m}^\star (\bx)]^T \bphi_{s,m+j}(\bx) \\
&\quad = \frac{1}{2} \int_{\RR^2}  \hspace{-0.05in}d \bx \, \left[ e^{m \tau \cL} \bphi_{s,0}^\star (\bx)\right]^T \bphi_{s,m+j}(\bx)\\
&\quad = \frac{1}{2} \int_{\RR^2}  \hspace{-0.05in}d \bx \, \left[  \bphi_{s,0}^\star (\bx)\right]^T e^{-m \tau \cL}\bphi_{s,m+j}(\bx) \\
&\quad = \frac{1}{2} \int_{\RR^2}  \hspace{-0.05in}d \bx \, \left[  \bphi_{s,0}^\star (\bx)\right]^T \bphi_{s,{2m+}j}(\bx) \\
&\quad \stackrel{\eqref{eq:A2}}{=} \frac{1}{2} \mathcal{D}_s(2 T+ (2m+j)\tau),
\end{align*}
and the proof of Theorem \ref{prop.2} is complete. $\Box$

%It remains to describe briefly the Gram-Schmidt orthogonalization procedure. 
%If we denote by $r_{i,j}$ the entries of $\bR_s$ and start with the first column in equation~\eqref{eq:I14},  we have 
%\[
%r_{0,0} v_{s,0}(\bx) = u_{s,0}(\bx), \quad 
%r_{0,0} = \|u_{s,0}\| = \sqrt{\lb u_{s,0},u_{s,0}\rb}.
%\]
%For the $j^{\rm th}$ column, with $j = 1, \ldots, M-1$, we have 
%\[
%\sum_{l=0}^j r_{l,j} v_{s,l}(\bx) = u_{s,j}(\bx), 
%\]
%where 
%\[
%r_{l,j} = \lb u_{s,j},v_{s,l} \rb, \quad 0 \le l \le j-1, \quad r_{j,j} = \| \mathbb{P}_{j-1} u_{s,j}\|,
%\]
%with $\mathbb{P}_{j-1}$ the orthogonal projection on the space spanned 
%by $v_{s,0}(\bx), \ldots, v_{s,j-1}(\bx)$.
%
%The Gram-Schmidt scheme gives an orthonormal basis
%\[
%\lb v_{s,j}, v_{s,l} \rb = \delta_{j,l}, \qquad \forall ~0 \le j, l \le M-1,
%\]
%that is causal, in the sense that for all $0 \le j \le M-1$ we have 
%\[
% v_{s,j}(\bx) \in \mbox{span} \{u_{s,0}(\bx), \ldots, u_{s,j}(\bx)\}.
% \]
\section{}
\label{ap:C}
To prove Theorem \ref{thm.2}, it is useful to write the solution of~\eqref{eq:Ev0}-\eqref{eq:Ev01} using functional calculus on the positive definite and self-adjoint operator $A(\hat{\eps})$,
\begin{equation}
u_s^{\rm ev}(\bx,t;\hat{\eps}) = \cos \big[ t \sqrt{A(\hat{\eps})} \big] u_{s,0}(\bx).
\label{eq:C1}
\end{equation}
The entries of the Gramian  are
\begin{align*}
\Big[\mathbb{G}_s^{\rm ev}(\hat{\eps})\Big]_{m,l} &= \lb  \cos \big[ m \tau  \sqrt{A(\hat{\eps})} \big] u_{s,0},\cos \big[ l \tau  \sqrt{A(\hat{\eps})} \big] u_{s,0}\rb \\
&= \lb u_{s,0}, \cos\big[ m \tau  \sqrt{A(\hat{\eps})} \big] \cos \big[ l \tau  \sqrt{A(\hat{\eps})} \big]u_{s,0}\rb 
\end{align*}
and due to the trigonometric identity of the cosine we have 
\begin{align*}
\Big[\mathbb{G}_s^{\rm ev}(\hat{\eps})\Big]_{m,l}&= \frac{1}{2} \lb u_{s,0}, \cos[(m+l)\tau \sqrt{A(\hat{\eps})} \big] u_{s,0}\rb \\ 
&+\frac{1}{2} \lb u_{s,0}, \cos[(l-m)\tau \sqrt{A(\hat{\eps})} \big] u_{s,0}\rb.
\end{align*}
The first term in the right hand side is, by equations~\eqref{eq:S2_C} and~\eqref{eq:C1}, the same as $\mathbb{D}_s^{\rm ev}( (m+l) \tau;\hat{\eps})$ 
and the second term is, similarly, $\mathbb{D}_s^{\rm ev}( (l-m) \tau;\hat{\eps})$. The result~\eqref{eq:Ev4} follows by setting $l = m +j$. 

To derive the identity~\eqref{eq:Ev5}, we observe from definitions~\eqref{eq:Sec0},~\eqref{eq:Ev1} and~\eqref{eq:S2_C} that 
\begin{align*}
\mathbb{D}_s^{\rm ev}(j \tau;\hat{\eps}) = \frac{1}{2}\left[\mathcal{C}_j(\hat{\eps}) + \mathcal{C}_{-j}(\hat{\eps})\right].
\end{align*}
Furthermore, equation~\eqref{eq:Sec1} gives 
\begin{align*}
\mathbb{D}_s^{\rm ev}(j \tau;\hat{\eps}) &= \frac{1}{2} \left[ D_s(j \tau) +  D_s(-j \tau)\right] \\
&+ \frac{1}{4} D_s(2T + j \tau;\hat{\eps}) + 
\frac{1}{4} D_s(2T - j \tau),
\end{align*} 
where we dropped the argument $\hat{\eps}$ in the last term because the data~\eqref{eq:S6}  evaluated at $t < 2T$ does not contain 
information about the heterogeneous part of the medium. We have 
\begin{align*}
\frac{1}{2} D_s(2T + j \tau;\hat{\eps}) &= 2 \mathbb{D}_s^{\rm ev}(j \tau;\hat{\eps}) - \frac{1}{2} D_s(2T - j \tau) \\
&- [D_s(j \tau) +  D_s(-j \tau)], \quad j \ge 0,
\end{align*}
and using this in the expression of $\mathbb{G}_s(\hat{\eps})$, given by the analogue of~\eqref{eq:I11} with $\eps(\bx)$ replaced by $\hat{\eps}(\bx)$, we get 
\begin{align*}
[\mathbb{G}_s(\hat{\eps})]_{m,m+j} = 2 \mathbb{D}_s^{\rm ev}((2m+j) \tau;\hat{\eps}) + (\Lambda_s)_{m,m+j},
\end{align*}
and 
\begin{align*}
[\mathbb{G}_s(\hat{\eps})]_{0,j} = 2 \mathbb{D}_s^{\rm ev}(j \tau;\hat{\eps}) + (\Lambda_s)_{0,j},
\end{align*}
with $\Lambda_s$ defined in~\eqref{eq:Ev6}.
The result~\eqref{eq:Ev5} follows from these identities and equation~\eqref{eq:Ev4}. $ ~ \Box$

\section{Details of numerical Implementation}
\label{ap:B}

%forward model
We solved the wave equation \eqref{eq:Ev0}-\eqref{eq:Ev01} with a second-order finite-difference time-domain (FDTD) method, truncated by an uniaxial perfectly matched layer (UPML). Following the anisotropic PML formulation of Gedney~\cite{TafloveHagness,Gedney1996}, we have a quadratically increasing delay profile and a quartically increasing decay profile. The maximum values of the delay and decay profiles were tuned to the pulse bandwidth and the expected angles of incidence. 
The domain $\Omega$ is discretized  with 25 points per wavelength $\la_o$ and we use  an explicit leap-frog integrator on a GPU for time stepping at intervals satisfying the CFL condition. 
%Adjoint computation

The minimization of the objective functions $\mathcal{O}( \mathcal{E}_{\balpha})$ and $\mathcal{O}^{\rm FWI}( \mathcal{E}_{\balpha})$ is done with a Gauss-Newton iteration  and a line search in the update directions, as described in~\cite{HabashyAbubakar2004}. The Jacobian is computed efficiently  using the well-known adjoint-state method \cite{Plessix2006}. This computes the sensitivities of the data with respect to  the components of $\balpha$, that parametrize the search permitivity \eqref{eq:Par1}, as a convolution in time of the wavefield $u_s(\bx,T+t)$ and its adjoint, and an inner product in space, over the support of the perturbation. This can be done in parallel for each slow time $s$. For our configurations, the wavefields of a single slow time $s$, at points in $\Omega_{\rm im}$,  fit entirely into GPU-VRAM memory.

\vspace{-0.1in} 
\subsection{Regularization} 
Our inversion approach requires taking the Cholesky square root of the data driven Gramian.  For noiseless data, this is a
symmetric and positive definite matrix by definition.  For noisy data, we replace $\mathbb{D}_s(0)$ by  $(1+\eta) \mathbb{D}_s(0)$
 in equation \eqref{eq:Ev4}.  Recall from \eqref{eq:S2_C} that $\mathbb{D}_s(0)$ is the norm squared of $u_{s,0}(\bx)$,
which is computed analytically, using formula \eqref{eq:N4} evaluated at $t = T$. Equation \eqref{eq:Ev4} gives a Gramian with  correct, symmetric Topelitz plus Hankel algebraic structure. The boosted initial data enhances its diagonal and, by  varying $\eta$, we can increase its smallest eigenvalue to a positive level.

Regularization is also needed to stabilize the Gauss-Newton iteration. We use a multiplicative regularization approach, following \cite{multInversDruskin,vandenBerg2001}:

\emph{Start with the initial guess $\balpha^{(1)}= {\bf 0}$. For $i \ge 1$, let $\balpha^{(i)}$ be the current iterate of the vector of coefficients. Solve 
\[
\min_{\delta \balpha} \mathcal{O}^{\rm lin}_{\balpha^{(i)}}(\delta \balpha) \mathcal{R}_{\balpha^{(i)}}(\delta \balpha) 
\]
where 
\[
\mathcal{O}^{\rm lin}_{\balpha^{(i)}}(\delta \balpha) \approx \mathcal{O}\Big(\mathcal{E}_{\balpha^{(i)} + \delta \balpha}\Big)
\]
approximates~\eqref{eq:Ev12} by  linearizing $\bR_s^{\rm ev}(\mathcal{E}_{\balpha + \delta \balpha}) (\bR_s^{\rm ev})^{-1}$ with respect to $\delta \balpha$. The regularization factor is 
\begin{align*}
 \mathcal{R}_{\balpha^{(i)}}(\delta \balpha) &= \int_{\Omega_{\rm im}} \hspace{-0.1in}d \bx \, \eta_{\balpha^{(i)}}(\bx) 
 \Big[ \Big| \sum_{q=1}^Q (\alpha_q^{(i)} + \delta \alpha_q)\nabla_{\bx} e^{- \frac{\|\bx-\bz_q\|^2}{2 \sigma^2}}\Big|^2 \\
 &\hspace{0.2in}+ \frac{1}{h^2} \mathcal{O}(\mathcal{E}_{\balpha^{(i)}})\Big], 
 \end{align*}
 where 
 \begin{align*}
 {
 \eta_{\balpha^{(i)}}(\bx) =
\frac{1}{|\Omega_{\rm im}|}
 \left[\Big| \sum_{q=1}^Q \alpha_q^{(i)} \nabla_{\bx} e^{- \frac{\|\bx-\bz_q\|^2}{2 \sigma^2}}\Big|^2  + \frac{1}{h^2}\mathcal{O}(\mathcal{E}_{\balpha^{(i)}})\right]^{-1}
 %.
 }
 \end{align*}
% \begin{align*}
% \eta_{\balpha^{(i)}}(\bx) =
% \left[\Big| \sum_{q=1}^Q \alpha_q^{(i)} \nabla_{\bx} e^{- \frac{\|\bx-\bz_q\|^2}{2 \sigma^2}}\Big|^2 |  + \frac{\mathcal{O}(\hat \eps_i)}{h^2}\right]^{-1}.
% \end{align*}
{Update the iterate as $\balpha^{(i+1)} = \balpha^{(i)} +l \delta \balpha,$ where l is determined by line search of the nonlinear objective $\mathcal{O}(\mathcal{E}_{\balpha^{(i)} + \delta \balpha})\mathcal{R}_{\balpha^{(i)}}(\delta \balpha) $.} Stop if the convergence criterium is reached.
}

Note that the only parameter that has to be chosen is $h^2$, the mesh size used in the parametrization of the search permittivity.  Increasing $h$ will increase the smoothness of the reconstruction. The multiplicative regularization term resets after each iteration and is unity for 
{$\delta \balpha ={\bf 0}$.}
%$\delta \alpha_q = \mathbf{0}$.

\vspace{-0.1in} 
\subsection{Computational cost} 
The cost  is dominated by the computation of the Jacobian of the objective functions. The adjoint state method \cite{Plessix2006} involves a temporal convolution integral of the fields  $u_s(\bx,T+t)$, which are computed on a fine time grid using FDTD on a GPU. The convolution is done in the frequency domain, via FFT.  For each coefficient in 
$\balpha^{(i)}$ we need to compute an integral of the time-convolved fields with the corresponding basis function. Then, the result is converted to the time domain via FFT. Once we have the sensitivity of the data to the components of $\balpha^{(i)}$, the mapping to the sensitivity of $\bR_s$ is done efficiently, via known formulas for the perturbation of the Cholesky square root.  The Jacobian calculation is  parallelizable over multi-GPU systems, since it is entirely independent for each slow time $s$. In our prototype implementation, we used a single GPU, so the computation time of the Jacobian was limited by the available VRAM memory-bandwidth of the used GPU. The computation time is further proportional to the acquisition length $M$, and the dimensionality of the problem after discretization with finite differences. The latter  scales quadratically with the domain size in 2D and cubically in 3D and becomes limiting when scaling to larger targets.

The computational domain used in our study is small when compared to realistic spot-SAR configurations, which can have $10^4-10^5$ center wavelengths in each spatial dimension. When scaling to realistic targets, iterating the objective function is unrealistic using our current hardware capabilities, since the Jacobian has to be recomputed after each iteration. However, for the first iteration, the Jacobian can be precomputed, and our method can be used as an imaging function similar to current SAR imaging strategies. 

\vspace{-0.1in}
\subsection{Noise model}\label{ap:Noise}
The noisy data model is 
\[
\mathbb{D}_s^{\rm n}(t)= \int_{\mathbb{R}^2} d \bx \int_{-T_b}^{T_b} dt'  J_s(\bx,-t') [u_s(\bx,t-t')+W(\bx,t-t')],
\]
where the random field $W(\bx,t)$ models additive noise. We choose a zero mean, stationary Gaussian process with  covariance function that is separable in space and time
\[
\mathbb{E}[ W(\bx + \Delta \bx, t + \Delta t) W(\bx ,t)] = \sigma_{{W}}^2 e^{- \frac{(\Delta \bx)^2}{\ell_x^2} - \frac{(\Delta t)^2}{\ell_t^2} }.
\]
We choose the correlation length $\ell_{\bx}=\lambda_o/2$, corresponding to typical sensor separation in phased arrays, and  the correlation time $\ell_t=0.7\, T_{\rm Ny}$, where $T_{\rm Ny}$ is the Nyquist period of the center frequency of the used pulse. The parameter $\sigma_{{W}}$ adjusts the noise level.

Due to the  averaging involved in the computation of $D_s(t)$, the net noise term in $\mathbb{D}_s^{\rm n}(t)$ is correlated {in time $t$}, but significantly weaker. This is precisely the reason why in SAR one uses long and chirped signals. 
%Even for the shorter signal used in our case, in the simulation results reported in Fig. \ref{fig:Cracks}, 
%the signal-to-noise ratio was increased from 0.3 to 20.
{In the simulation results reported in Fig. \ref{fig:Cracks}, in which the source pulse is short, 
the signal-to-noise ratio for the signal $u_s+W$ in the support of $J_s$ is $0.3$, which gives a signal-to-noise ratio for the noisy data $\mathbb{D}_s^{\rm n}(t)$ equal to $20$.}

% --------------
% this is not a bib file so I have no idea why there is a IEEETRAN bibliography style here. it does not do anything.
\bibliographystyle{IEEEtran}
%Lilianas Bib:
%\include{bibliography}
% Jorns bib: fixed author spelling and style to IEEE
{

}
\end{document}